\RequirePackage[l2tabu, orthodox]{nag}
\documentclass[a4paper]{amsart}

\usepackage[T1]{fontenc}
\usepackage{cfr-lm}
\usepackage{microtype}

\usepackage{amsmath}
\usepackage{amsthm}
\usepackage{amssymb}
\usepackage[foot]{amsaddr}
\usepackage{mathtools}

\usepackage{tikz}
\usetikzlibrary{arrows}
\usetikzlibrary{cd}
\usepackage{physics}
\usepackage{xcolor}
\usepackage{graphicx}

\usepackage{subcaption}

\DeclareCaptionLabelFormat{xx}{#2}
\usepackage{listings}
\usepackage{booktabs}
\usepackage{enumitem}
\usepackage[backend=bibtex,style=numeric,autocite=inline, maxnames=8, isbn=false, sortcites=true, backref]{biblatex}
\usepackage{hyperref}
\hypersetup{
    colorlinks,
    linkcolor={red!50!black},
    citecolor={blue!50!black},
    urlcolor={blue!80!black}
}

\usepackage[capitalise,noabbrev]{cleveref} 

\newtheorem{thm}{Theorem}
\numberwithin{thm}{section}
\newtheorem*{thm*}{Theorem}

\AddToHook{env/lem/begin}{\crefalias{thm}{lem}}

\AddToHook{env/prp/begin}{\crefalias{thm}{prp}}

\AddToHook{env/cor/begin}{\crefalias{thm}{cor}}

\AddToHook{env/obs/begin}{\crefalias{thm}{obs}}
\theoremstyle{definition}
\newtheorem{defn}[thm]{Definition}
\AddToHook{env/defn/begin}{\crefalias{thm}{defn}}
\newtheorem{ex}[thm]{Example}
\AddToHook{env/ex/begin}{\crefalias{thm}{ex}}
\newtheorem{problem}[thm]{Problem}
\AddToHook{env/problem/begin}{\crefalias{thm}{problem}}
\newtheorem{cons}[thm]{Construction}
\AddToHook{env/cons/begin}{\crefalias{thm}{cons}}
\newtheorem{alg}[thm]{Algorithm}
\AddToHook{env/alg/begin}{\crefalias{thm}{alg}}
\theoremstyle{remark}
\newtheorem{rem}[thm]{Remark}
\AddToHook{env/rem/begin}{\crefalias{thm}{rem}}

\newcommand{\df}{\textit}

\newcommand{\union}{\cup}
\newcommand{\inter}{\cap}

\newcommand{\mc}{\mathcal}

\newcommand{\mf}{\mathfrak}

\newcommand{\git}{\mathord{
  \mathchoice{/\mkern-6mu/}
    {/\mkern-6mu/}
    {/\mkern-5mu/}
    {/\mkern-5mu/}}}

\renewcommand{\H}{\mathbb{H}}

\newcommand{\R}{\mathbb{R}}
\newcommand{\C}{\mathbb{C}}

\newcommand{\Z}{\mathbb{Z}}
\newcommand{\Q}{\mathbb{Q}}

\newcommand{\Sph}{\mathbb{S}}

\DeclareMathOperator{\PSL}{\mathsf{PSL}}

\DeclareMathOperator{\Isom}{Isom}

\DeclareMathOperator{\Hom}{Hom}

\DeclareMathOperator{\Arg}{arg}

\makeatletter
\@namedef{subjclassname@2020}{%
  \textup{2020} Mathematics Subject Classification}
\makeatother

\title[From disc patterns to knot groups]{From disc patterns in the plane to character varieties of knot groups}
\author[A. Elzenaar]{Alex Elzenaar}
\address{School of Mathematics, Monash University, Melbourne, Australia}
\email{alexander.elzenaar@monash.edu}
\thanks{\noindent This research was supported in part by the Monash eResearch Centre and eSolutions-Research Support Services through the use of the MonARCH HPC Cluster. The
author was supported by an Australian Government Research Training Program (RTP) Scholarship during part of the period that this work was undertaken. I thank the referee for
helpful comments.}

\keywords{computational geometry, Kleinian groups, discrete groups, hyperbolic $3$-manifolds, Schottky groups, unknotting tunnels, fractals, quasiconformal deformation spaces}
\subjclass[2020]{Primary 57K31; Secondary 20H10, 30F40, 52C25, 52C26, 57-08, 57K10, 57K32, 57K35, 58H15}

\begin{document}
\begin{abstract}
Motivated by an experimental study of groups generated by reflections in planar patterns of tangent circles, we describe some methods for constructing
and studying representation spaces of holonomy groups of infinite volume hyperbolic $3$-manifolds that arise from unknotting tunnels of links. We include full descriptions of
our computational methods, which were guided by simplicity and generality rather than by being particularly efficient in special cases. This makes them easy for non-experts to
understand and implement to produce visualisations that can suggest conjectures and support algebraic calculations in the character variety. Throughout, we have tried to
make the exposition clear and understandable for graduate students in geometric topology and related fields.
\end{abstract}

\maketitle

\section*{Introduction}
This article is a potpourri of different constructions of families of Kleinian groups, with an emphasis on computer experiment and visualisation. Instead of proving general
and technical theorems for an audience of experts, we have tried to give examples that are interesting to graduate students and to researchers in related fields (knot theory,
computational geometry, low-dimensional topology). We were influenced by the pictures and examples in the classic works of Fricke and Klein \autocite{fricke_klein65,fricke_klein66},
the work of Riley in the 1970s on knot group representations \autocite{riley72,riley79,riley75,riley75c,riley75b,Riley13} (see also the bibliography of his work appended
to \autocite{Brin13}), and the work of Sakuma and various co-authors on unknotting tunnels \autocite{yokota96,sakuma98,akiyoshi,lee12}. Pulling together some common threads
from all of these sources, we will travel in a straight line (with some minor detours) from the classical study of discrete groups generated by reflections in symmetric patterns of circles
to a modern view of some aspects of the representation theory of knot complement holonomy groups. Some of our work that we touch on here was motivated by the writing
of Hodgson and Kerckhoff on deformations of cone manifold structures \autocite{hk98,hk03}, and in particular questions about whether such deformations can be achieved via entirely
classical geometric methods, but our work in this direction is only at an early stage.

In \cref{sec:necklace}, we study groups arising from closed loops of tangent discs in the plane (necklaces), and we give an explicit construction of
an infinitely generated discrete group that can be implemented on the computer. Readers primarily interested in the knot theory may safely skip the
details of the two labelled constructions in this section, and skim it only to pick up language and notation. We then proceed, in \cref{sec:lattices}, to study groups which come from
lattice disc packings of the plane and we observe interesting connections to the theory of knots. This motivates our study in \cref{sec:families} of slices
through character varieties, and in \cref{sec:paths} of non-linear parameterised paths that lead from knot groups to the boundary of certain Teichm\"uller spaces.
Some aspects of this last study were discussed in the author's talk at the Hodgsonfest.

Most of the computer code used to produce the images in this article is available online at \url{https://github.com/aelzenaar/disc-patterns-and-knots}.

\section{Necklaces}\label{sec:necklace}
The history of groups generated by circle reflections in the plane is a long and interesting one; for a detailed discussion of the origins of the theory see the excellent
history of complex analysis by Bottazzini and Gray \autocite{bottazzini_gray}, Gray's biography of Poincar\'e \autocite{grayP}, and Stillwell's annotated translations
of original papers by Poincar\'e and Dehn \autocite{poincare,dehn}. In the late 19th century, a pair of influential books on the fledgling subject was written
by Klein and Fricke \autocite{fricke_klein65,fricke_klein66} (translated to English as \autocite{fricke_klein65t,fricke_klein66t}), including many pictures of reflection groups
acting on $ \Sph^2 $, $ \H^2 $, and $ \R^2 $ which have inspired generations of mathematicians (see \autocite{indras_pearls}, the preface to \autocite{magnus}, and part III
of \autocite{mandelbrotCW}), researchers in adjacent fields (for examples from chemistry, see \autocite{hyde03}), and artists (see e.g.\ \autocite{stier17,rigby95,bulatov13}).
We will discuss some aspects of one of their pictures, \autocite[Vol. I, fig. 156]{fricke_klein65} (reproduced as \autocite[Fig.~65 on p.~133]{kag}), which shows a jagged
topological curve joining up the tangency points of a chain of tangent circles. Apart from the aesthetic beauty of the groups we will construct, this study will allow us
to introduce some of the basic definitions which we use throughout the paper. As an application, we will give some computationally realisable approximations to an interesting
infinitely generated discrete group.

\begin{defn}\label{defn:reflection}
  Let $ C $ be a circle in the plane, with centre $ w $ and radius $ r $. The \df{reflection} across $ C $ is the map which sends $ z $ to the point $ z' $ which
  lies on the ray $ \overrightarrow{wz} $ such that $ \abs{w - z} \abs{w - z'} = r^2 $;  explicitly, it is the anti-meromorphic function
  \begin{displaymath}
    z \mapsto \frac{r^2}{\bar{z} - \bar{w}} + w.
  \end{displaymath}
  A cyclically ordered set of points $ w_1,\ldots,w_n \in \C $ determines a (possibly self-intersecting) polygon. Draw circles
  $ C_i $ centred at each $ w_i $, with radii chosen so that $ C_i $ is tangent to $ C_{i \pm 1} $ for all $ i $. Let $ \Delta_i $ be the open disc bounded by $ C_i $
  (i.e.\ the bounded component of $ \C \setminus C_i $), and suppose that the radii are chosen so that if $ \Delta_i $ and $ \Delta_j $ intersect for any $ i $ and $ j $,
  then $ \Delta_i = \Delta_j $ (so if discs overlap then they coincide). Write $ \phi_i $ for the reflection in the circle $ C_i $, and let $ \tilde{\Gamma} = \langle \phi_1,\ldots,\phi_n \rangle  $.
  The orientation-preserving half of $ \tilde{\Gamma} $, which we denote $ \Gamma $, is $ \langle \phi_i \phi_{i+1} : 1 \leq i \leq n  \rangle $ (subscripts taken mod $ n $).
  It is called a \df{necklace group} or a \df{bead group} \autocite[\S VIII.F]{maskit}; the `necklace' is the chain of beads $ \{C_1,\ldots,C_n\} $.
\end{defn}

\begin{figure}
  \begin{subfigure}{0.49\textwidth}
    \centering
    \includegraphics[width=\textwidth]{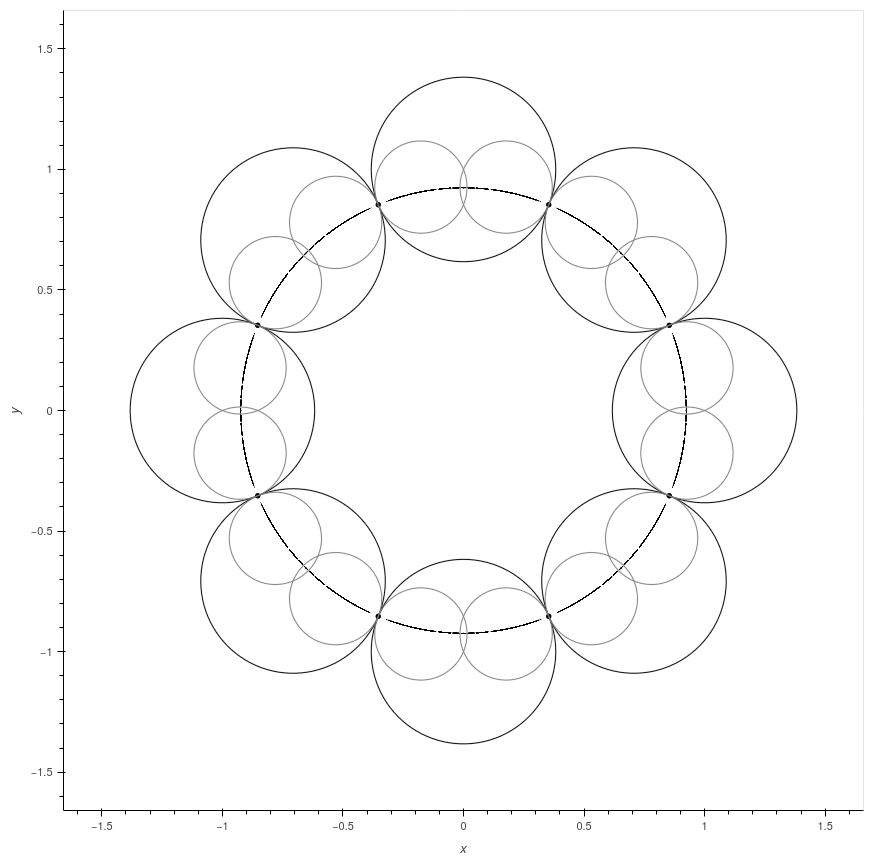}
    \caption{Roots of unity}
  \end{subfigure}\hfill%
  \begin{subfigure}{0.49\textwidth}
    \centering
    \includegraphics[width=\textwidth]{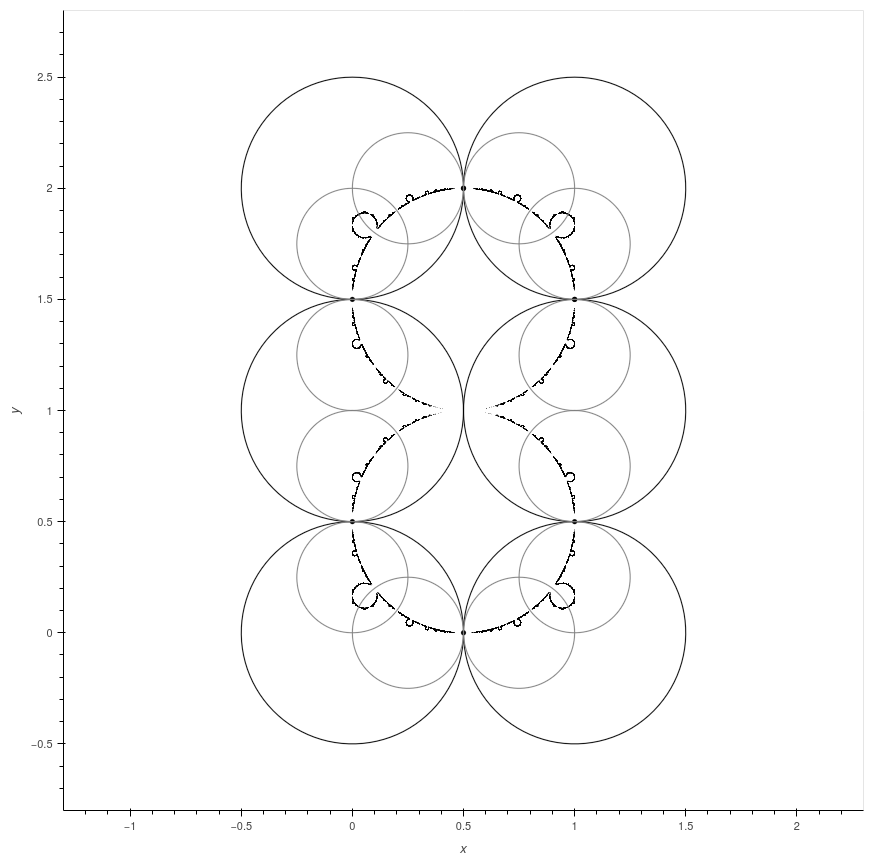}
    \caption{Rectangle}
  \end{subfigure}\\[2ex]
  \begin{subfigure}{0.49\textwidth}
    \centering
    \includegraphics[width=\textwidth]{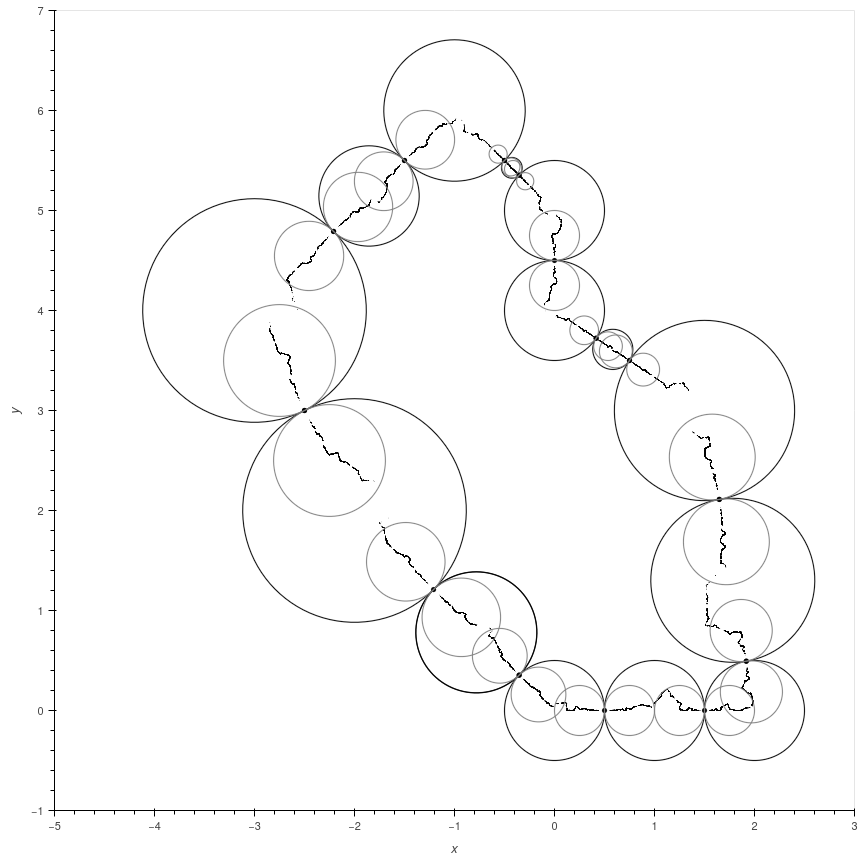}
    \caption{A complicated polygon}
  \end{subfigure}\hfill%
  \begin{subfigure}{0.49\textwidth}
    \centering
    \includegraphics[width=\textwidth]{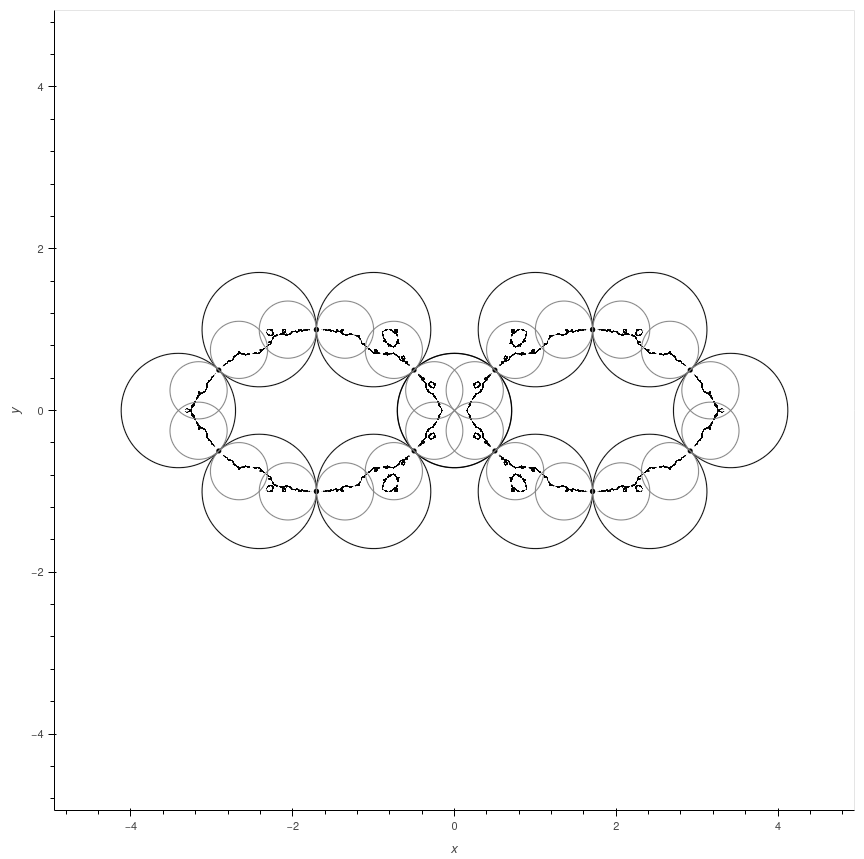}
    \caption{Figure 8}
  \end{subfigure}
  \caption{Four necklace groups.\label{fig:bead_groups}}
\end{figure}

It follows from the Poincar\'e polyhedron theorem that $ \Gamma $ is a discrete group of conformal automorphisms of the Riemann sphere.
We will give several explicit examples in \cref{ex:necklaces} below, and the definition we just gave might be clarified by looking at \cref{fig:bead_groups}. The necklace groups are our
first family of examples of \df{Kleinian groups}, which are equivalently (i) discrete groups of orientation-preserving conformal automorphisms of $ \Sph^2 $; (ii) discrete
subgroups of $ \PSL(2,\C) $, acting as fractional linear transformations on the Riemann sphere $ \hat{\C} = \C \union \{\infty\} $; (iii) discrete subgroups of $ \Isom^+(\H^3) $.
The elementary theory of Kleinian groups may be found in Maskit \autocite{maskit}, and any undefined and unfamiliar words may be looked up there.

Our examples use the following explicit construction of a necklace group from a polygon. It introduces additional
vertices if needed so that the circles in the necklace are disjoint.
\begin{cons}[Finite necklace groups]\label{cons:necklace}
  The construction takes as input a tuple $ (p_0,\ldots,p_n) $ of complex numbers and outputs a list $ \mc{L} $ of circles.
  \begin{enumerate}
    \item Refine the polygon defined by the $ (p_i) $ by adding vertices at any self-intersection points, and then adding interpolating vertices in the
          interiors of the edges so that the closest vertices to a given one are exactly those which are joined to it along a polygon edge. Replace the tuple $ (p_i) $
          with this larger ordered set of points.
    \item For each $ p_i $, set $ r_i = \frac{1}{2}\min \{ \abs{p_{i-1}-p_i}, \abs{p_i - p_{i+1}} \} $. Add to $ \mc{L} $ the circle $ C_i $ centred at $ p_i $ of radius $ r_i $.
    \item For each $ i $ consider the circles $C_i$ and $C_{i+1}$. By construction, these circles do not intersect, except possibly at a point of tangency.
          It might be that the circles are not tangent, and we add an additional circle to correct this. Let $x$ and $y$ be the points of the segment $[p_i, p_{i+1}]$
          which are respectively the intersection points of $C_i$ and $C_{i+1}$ with that segment. Append to $ \mc{L} $ the circle centred at $(x+y)/2$ with radius $\abs{x-y}/2$.
  \end{enumerate}
  The group generated by the reflections in the circles contained in $ \mc{L} $ is discrete.
\end{cons}

Using the construction, one can obtain interesting $\R$-algebraic maps from configuration spaces of rigid bars and joints (which are of great interest in real algebraic
geometry and kinematics, see Cox \autocite[Chapter~4]{coxAPS}) to deformation spaces of Kleinian groups. A very simple map of this type is described in \cref{ex:pendulum}.

\begin{ex}[Four necklace groups]\label{ex:necklaces}
  In \cref{fig:bead_groups}, we show (in black) the defining circles $ C_i $ for four different necklace groups; the points chosen in each case are respectively:
  \begin{enumerate}[label=(\Alph*)]
    \item The $8$th roots of unity;
    \item $\{0, 1, 1+1i, 1+2i, 2i, 1i\} $;
    \item $ \{0, 1, 2, 1.5+3i, 4i, 5i, 6i - 1, 4i - 3, 2i-2\}$; and
    \item $\{0, 1+1i, 1+\sqrt{2} + 1i, 2+\sqrt{2}, 1 + \sqrt{2} - 1i, 1-1i, 0, -1+1i, -1-\sqrt{2} + 1i, -2-\sqrt{2}, -1 - \sqrt{2} - 1i, -1-1i\}$ (compare Maskit \autocite[\S VIII.F.5]{maskit}).
  \end{enumerate}

  The light gray circles are the isometric circles of the chosen generators for the orientation-preserving half. At each point of tangency of the black reflection circles there is also
  a tangent pair of isometric circles, and the product of the reflections in the two black circles maps the interior of one grey circle onto the exterior of the other grey circle.

  The jagged black curves are the \df{limit sets} of the necklace groups, which are fractals analogous to Julia sets from complex function
  dynamics: if $ \Gamma $ is a group of conformal automorphisms of $ \hat{\C} $, then its limit set is the set $ \Lambda(\Gamma) $ of all accumulation points of orbits of points of
  $ \hat{\C} $ under $ \Gamma $. By taking sufficiently small vertex circles and sufficiently small interpolating circles in the construction, one can obtain groups whose limit sets
  are arbitrarily close to the starting polygon.
\end{ex}

The necklace groups and their limit sets can be defined for groups of conformal automorphisms acting on higher dimensional spheres. For instance, in $ \Sph^3 $ taking
the group generated by reflections in tangent spheres along a knotted polygon gives a limit set that is a knotted fractal curve (in fact, a wild knot); this construction
is originally due to Apanasov \autocite{apanasov80r,apanasov80e} and a physical model appears in Díaz, Hinojosa, Mendoza, and Verjovsky \autocite{diaz19}.
One can even construct groups in $ \Sph^3 $ whose limit sets are wild knotted $2$-spheres, see Apanasov \autocite[Example~7.19]{apanasov}.

Methods for drawing the limit set of a Kleinian group $\Gamma$ tend to rely on the fact that the action of $ \Gamma $ on its limit set is ergodic.
An extensive and readable study of the computational aspects of limit sets is that by Mumford, Series, and Wright \autocite{indras_pearls}. Further references on
fractals associated to Kleinian groups may be found in Nakamura \autocite{nakamura21}.

The complement of the limit set of a group $ \Gamma \leq \PSL(2,\C) $ in $ \hat{\C} $ is called the \df{domain of discontinuity} and will be denoted by $ \Omega(\Gamma) $. It is
analogous to the Fatou set of a dynamical system. If $ \Omega(\Gamma) \neq \emptyset $ then the group is Kleinian, but the converse is not true; for instance if $ \Gamma $ is the
holonomy group of a hyperbolic knot complement then $ \Lambda(\Gamma) = \hat{\C} $. The quotient $ \Omega(\Gamma)/\Gamma $ is a  (possibly singular) Riemann surface that can be
naturally identified with the conformal boundary of the orbifold $ \H^3/\Gamma $. We will often say that $ \Gamma $ \df{uniformises} $ \Omega(\Gamma)/\Gamma $; this is an old usage
of the word, which arises from the early work of Riemann and Poincar\'e on automorphic functions---$ \Gamma $ was viewed by them as a gadget by which functions defined on a surface
could be lifted to a planar set, and this general process was called uniformisation of the surface.

Our own motivation for implementing a general procedure to get necklace groups from polygons was to produce visualisations of some remarkable
groups first studied by Accola.
\begin{ex}[Atom groups]
  In 1966, Accola \autocite{accola66} constructed a Kleinian group $A$ with $ \Omega(A) \neq \emptyset $ that is
  provably not finitely generated. Let $ \mc{A} $ be the annulus bounded by the circles of radius $1$ and $2$ centred at $ 0 $.
  Draw two disjoint spirals $ S_1, S_2 $ in $ \mc{A} $ which limit at each end onto the boundary circles of $ \mc{A} $. For each spiral $ S_i $ ($i \in \{1,2\} $)
  draw a sequence of mutually tangent circles which (i) bound discs whose union covers $ S_i $ and which lies in $ \mc{A} $, and which (ii) do not intersect any of the circles chosen on $ S_{3-i} $.
  Now let $ A $ be the group generated by the reflections in these circles. The corresponding quotient surface has four components, including two discs (the projections of the exterior discs of $ \mc{A} $),
  and by Ahfors' finiteness theorem \autocite[\S 4.2]{matsuzakitaniguchi} a finitely generated group cannot uniformise discs.
\end{ex}

\begin{figure}
  \begin{subfigure}[t]{0.49\textwidth}
    \centering
    \includegraphics[width=\textwidth]{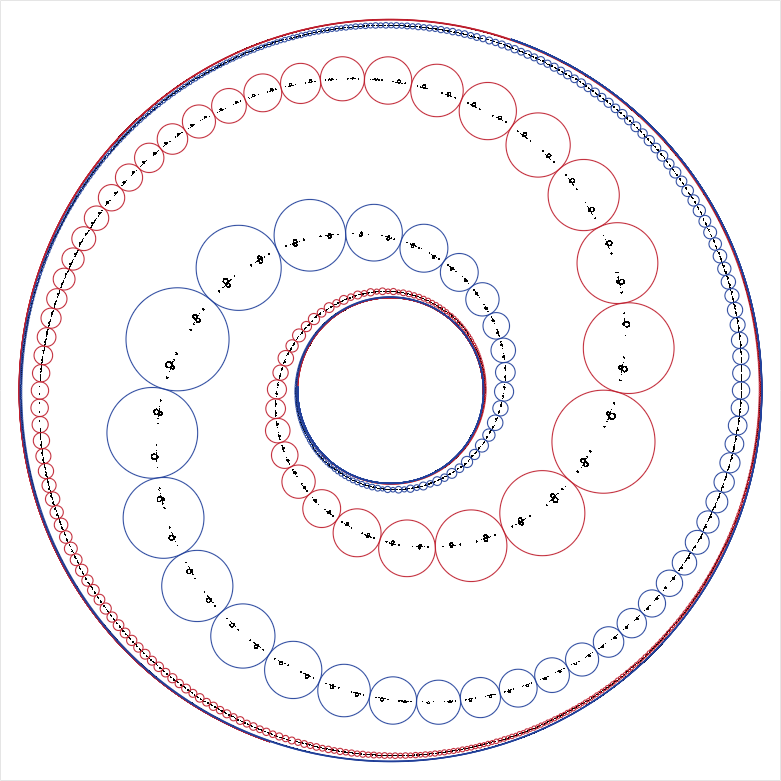}
    \caption{}
  \end{subfigure}\hfill%
  \begin{subfigure}[t]{0.49\textwidth}
    \centering
    \includegraphics[width=\textwidth]{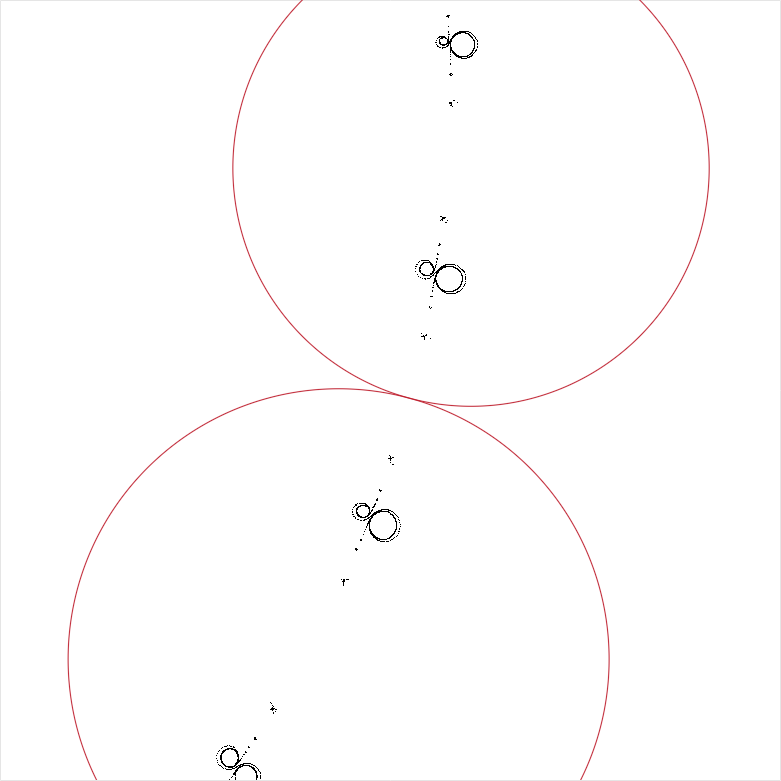}
    \caption{}
  \end{subfigure}
  \caption{An approximation to an Accola atom group. The group is generated by the reflections in the circles arranged along the two spirals. The limit set
           is approximated by the points visible within the circles. If the picture had infinite resolution, we would be able to see that the limit set is
           itself made up recursively of spirals.\label{fig:limsets3}}
\end{figure}

Realising approximations to this group on the computer takes a bit of care. We approximate the two spirals locally by circles and use trigonometry to place vertices for
the application of \cref{cons:necklace}. The problem is that, in practice, one ends up with a picture where massive circles alternate with tiny circles. We solve this by `nudging'
the radii of the circles. The result of the construction with 10,000 generators is shown in \cref{fig:limsets3}.
\begin{cons}[An explicit atom group]
  Let $ \sigma_1(\theta) $ and $ \sigma_2(\theta) $ be the polar parameterisations of the two spirals. For example,
  \begin{displaymath}
    \sigma_1(\theta) = \frac{1}{2} + \frac{3}{2+2\exp(-0.8\theta)} \;\text{and}\; \sigma_2(\theta) = \sigma_1(\pi + \theta).
  \end{displaymath}
  Let $ 2n $ be the number of circles to produce along each spiral in the approximation. For each $ i \in \{1,2\} $ run the following algorithm to place circles along the spiral $ \sigma_i $; our starting
  angle $ \theta_0 $ is defined to be $ \theta_0 = \pi $ for $ i = 1 $ and $ \theta_0 = 0 $ for $ i = 2 $.
  \begin{enumerate}
    \item Define the function
      \begin{multline*}
        r_{\text{guess}}(\theta) = \frac{1}{2.2} \min\Bigl\{\, \abs{\sigma_{i+1}(\theta - 2\pi) - \sigma_i(\theta)},\\\abs{\sigma_{i+1}(\theta+2\pi)- \sigma_i(\theta)}, \abs{\sigma_{i+1}(\theta) - \sigma_i(\theta)}\,\Bigr\},
      \end{multline*}
      which assigns to each angle $ \theta $ a radius which ensures that it doesn't collide with circles placed on other arms of the spiral. The exact factor $ 2.2 $ was chosen for aesthetic reasons.
    \item Set $ c = \sigma_i(\theta_0) \exp(\theta_0 i) $ and $ r = r_{\text{guess}}(\theta_0) $. Let $ \mc{C} $ be a double-ended queue, initialised with the single triple $ (\theta_0, r, c) $.
    \item Iterate the following $ n $ times with $ \delta = +1 $ and $ n $ times with $ \delta = -1 $:
      \begin{enumerate}
        \item If $ \delta = -1 $ (resp.\ $ +1$) then let $ \theta_{\text{prev}} $, $ c_{\text{prev}} $, and $ r_{\text{prev}} $ be the three elements of the first (resp.\ last) entry of $ \mc{C} $.
        \item Set $ r_\text{target} = r_{\text{guess}}(\theta_{\text{prev}}) $.
        \item Set $ \Delta \theta =  -\delta\arctan\left( \frac{r_\text{target}}{\sigma_i(\theta_{\text{prev}})}\right) $. This is the approximate change in angle one would see
              from one circle centre to the next at $ \theta $ if we were placing them around a circle and not a spiral, and if we took radii to be exactly $ r_{\text{guess}} $.
        \item Set $ l = 0 $ and $ r = 0 $. We will modify $ r $ so that it ends up close to $ r_\text{target} $ by nudging $ \Delta\theta $ up and down slightly.
        \item While $ \abs{r/r_{\text{target}}} < .8 $ or $ \abs{r/r_{\text{target}}} > 1.2 $:
          \begin{enumerate}
            \item If $ r < r_{\text{target}} $ then set $ \Delta\theta = 1.2l\Delta\theta $. Otherwise set $ \Delta\theta = .8 l\Delta\theta $.
            \item Set
              \begin{displaymath}
                c = \exp(i(\theta_{\text{prev}} + 2\Delta\theta) )\, \sigma_i(\theta_{\text{prev}} + 2\Delta\theta)
              \end{displaymath}
                and $r = \big\lvert\,\abs{c_{\text{prev}}-c} - r_\text{prev}\big\rvert$.
            \item Increment $ l $.
          \end{enumerate}
        \item If $ \delta = -1 $ (resp.\ $ +1$) then prepend (resp.\ append) to $ \mc{C} $ the triple $ \left(\theta_{\text{prev}} + 2\Delta\theta, c, r\right) $.
      \end{enumerate}
    \item The set of centres and radii stored in $ \mc{C} $ now give circles whose reflections approximate the Accola atom group.
  \end{enumerate}
\end{cons}

Initial experiments in trying to construct these groups on the computer teach an important lesson about how challenging it is to work with these
geometric constructions when the groups are given in terms of $ 2 \times 2 $ matrices: sometimes it is much easier to work directly with moduli
coming from the combinatorics and geometry of circles, and only perform the transformation into matrix algebra once these moduli are settled,
rather than trying to write all the moduli in terms of equations involving entries of matrices: in other words, to work with the algebraic geometry of the
action on the space of circles, rather than the algebraic geometry of the representations. We will see this in the next section, where we use global
symmetries of a fundamental domain to reduce the dimension of the representation spaces we consider.

\section{Circle packings and knots}\label{sec:lattices}
Now that we have warmed up by constructing some groups from the combinatorial data of circle chains, we will look at groups coming from systems of circles with more global symmetry.
This symmetry will imply that the character varieties of these groups are algebraically very simple. As we alluded to in the previous section, the literature of specific examples of discrete
groups in $ \PSL(2,\C) $ is incredibly rich. One particularly interesting vein of examples was studied by Wielenberg in the 1970s \autocite{wielenberg78} (see also \autocite[Examples~59--62]{kag}).
He was interested in constructing finite-index subgroups of the Picard group, and much of the work we will describe here was directly inspired by his examples.

It is well-known to small children that there are not very many ways to pack the plane with equally sized discs in a way which is invariant under the action of a lattice. The simplest Kleinian
groups which can be defined in terms of this combinatorial structure are generated by three elements: two translations which generate the lattice, and a single additional element which maps a
pair of neighbouring circles into each other. Generically, when the circles are small compared to the lattice, these groups uniformise $ (1;2)$-compression
bodies \autocite{elzenaar24c,dang19,lackenbypurcell13}. Up to conjugacy, they are determined by three complex parameters and their character variety is parameterised by the map
\begin{displaymath}
  \C^3 \ni (\alpha,\beta,\lambda) \mapsto G(\alpha,\beta,\lambda) \coloneq
  \left\langle\!
  P = \begin{bmatrix} 1 & \alpha \\ 0 & 1 \end{bmatrix}\!, \;
  Q = \begin{bmatrix} 1 & \beta \\ 0 & 1 \end{bmatrix}\!,\;
  M = \begin{bmatrix} \lambda & \lambda^2 - 1 \\ 1 & \lambda \end{bmatrix}
 \right\rangle\!.
\end{displaymath}

\begin{figure}
  \centering
  \includegraphics[width=.5\textwidth]{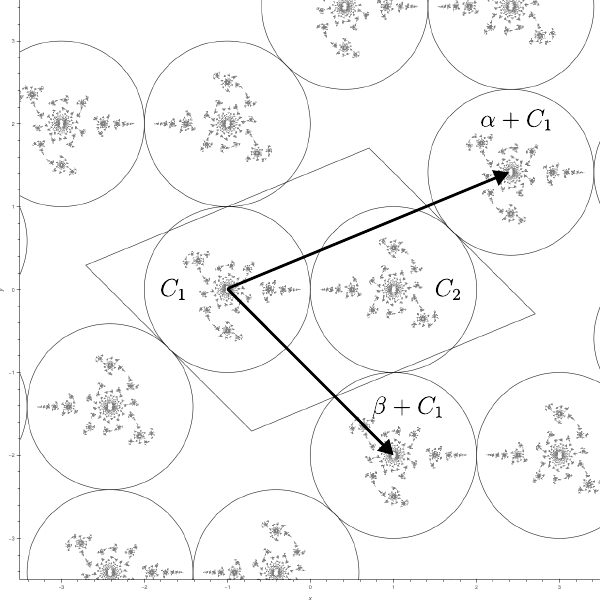}
  \caption{The element $ M $ and its conjugates do not preserve the set of tangencies of this asymmetric circle pattern. As a consequence, we see from the limit set (in grey) that many of the
          tangency points do not correspond to parabolic fixed points.\label{fig:bad_hexagon}}
\end{figure}

If we pick a random triple in $ \C^3 $ then there is no reason for the corresponding group to even be discrete, let alone torsion-free. We will place some simple
combinatorial conditions on the isometric circles of $ M $ to make them compatible with the lattice. Let $ C_1 $ and $ C_2 $ be the circles of radius $1$ centred at $-1 $ and $ +1 $
respectively, and let $ D $ be the union of the two circles and their interiors. Fix a lattice $ \Lambda \subset \C\oplus \C $ of rank $ 2 $, marked with
two generators $ \alpha $ and $ \beta $ such that
\begin{enumerate}
  \item[(CP1)] $ \alpha + C_1 $ and $ \beta + C_1 $ are tangent to $ C_2 $, and
  \item[(CP2)] if $ C_i \inter (\lambda + C_j) \neq \emptyset $ for any $ \lambda \in \Lambda $ and any $ i,j \in \{1,2\} $, then
               either $ C_i = \lambda + C_j $ or $ C_i $ is tangent to $ \lambda + C_j $.
\end{enumerate}
See \cref{fig:bad_hexagon} for an example illustrating these conditions. They imply that $ \Lambda + D $ is a disc pattern in $ \C $ (in the sense of Gr\"unbaum and Shephard \autocite[\S 7.3]{grunbaum}).
The complement of $ \Lambda + D $ is a union of topological polygons with circular arcs as edges and vertices coming from tangency points between circles. There are between three
and six circles (including $ C_1 $) tangent to $ C_2 $. The centres of all the circles in the pattern form the vertices of a hexagonal tiling of the plane of edge length $2$ (we use
the word `tiling' loosely, since it could be that the tiles overlap if the lattice translations are short).

We will consider groups of the form $ G(\alpha,\beta,\lambda) $ with the property that $ M $ is a parabolic which pairs the two circles $ C_1 $ and $ C_2 $: a sufficient condition
for this is that $ \lambda = 1 $, in whch case $ C_1 $ and $ C_2 $ are the isometric circles of $ M $. Denote by $ \mc{Q} $ the quadrilateral with vertices $ (\pm \alpha \pm \beta)/2 $
(where the two signs are independent); this is a fundamental domain for the action of the lattice $ \Lambda $. It is not necessarily the case that the common
exterior of all the circles, intersected with $ \mc{Q} $, gives a fundamental domain for the action of the group $ G(\alpha,\beta,\lambda) $. We need to impose further
conditions that force $ P $, $ Q $, and $ M $ to form a side-pairing structure on the polygons; and then we must verify that the vertex cycle elements are parabolic \autocite[\S II.G]{maskit}.

For the side-pairings to be well-defined, we need $ M $ to send the point of tangency of $ C_1 \inter (-\beta + C_2) $ to the point $ C_2  \inter (\alpha + C_1) $, and to
send the point of tangency $ C_1 \inter (-\alpha + C_2) $ to the point $ C_2  \inter (\beta + C_1) $ (\cref{fig:bad_hexagon} shows an example where this fails). Since $ M $
is parabolic, this is equivalent to imposing the following additional condition on the lattice:
\begin{enumerate}
  \item[(CP3)] $\mc{Q}$ must be symmetric with respect to reflection across $ i\R $. Hence $ \alpha = \pm\overline{\beta} $.
\end{enumerate}
There are two parabolic cycles in these side-pairing structures if the hexagon has not degenerated to form additional tangencies, corresponding to the elements $ M $ and $ M P^{-1} M Q^{-1} $.
The trace of the latter is $ (\alpha-2)(\beta-2) - 2 $; this trace is equal to $ -2 $ if and only if $ \alpha = 2 $ or $ \beta = 2 $, and since $ \alpha = \pm \overline{\beta} $
this would imply $ \alpha $ and $ \beta $ are linearly dependent. Thus a necessary condition for $M P^{-1} M Q^{-1}$ to be parabolic if the lattice $ \Lambda $ is rank $2$ is
that $ \tr M P^{-1} M Q^{-1} = 2 $, i.e.
\begin{equation}\tag{$\ast$}\label{eq:alpha}
  0 = (\alpha-2)(\beta-2) - 4 = \abs{\alpha}^2 - 4\Re(\alpha).
\end{equation}
The points of the lattice $ \Lambda = \Z \alpha + \Z \beta $ form the vertices of a number of hexagonal cells coming from cycles of six tangent circles. These cells are possibly overlapping, if
the values of $ \alpha $ and $ \beta $ are small (an example like this appears in \cref{fig:circles_2pi3}). Let $ \delta $ be the dihedral angle at the vertex $ +1 $ of the hexagon
with vertices $ \pm 1 $, $ -1 + \alpha $, $ 1-\beta $, $ 1-\beta+\alpha $, and $ -1+\alpha-\beta $. Then $ \alpha = 2 + 2\cos \delta + 2i\sin\delta $, so
\begin{displaymath}
  \abs{\alpha}^2 = (2 + 2\cos\delta)^2 + 4\sin^2\delta = 4 + 8\cos\delta + 4 = 4(2 + 2\cos\delta) = 4\Re(\alpha)
\end{displaymath}
and every value of $ \delta $ satisfies \eqref{eq:alpha}. In other words, the symmetry of $ \mc{Q} $ is also sufficient to imply that the two vertex cycles $ M $
and $ M P^{-1} M Q^{-1} $ are parabolic, which is enough to use the Poincar\'e polyhedron theorem to see that the group defined is discrete if the only tangent circles
are those that are coming from the cycles of six tangent circles around each hexagon. If the hexagons are degenerate (e.g.\ become rectangles) or the lattice is small (so they begin to overlap), then
additional points of tangency can arise. It turns out that the symmetry of $ \mc{Q} $ is sufficient to ensure that all these additional possible tangency points
are also parabolic fixed points, so no further conditions are needed. There are only finitely many values of $ \delta $ satisfying condition (CP2) which give rise
to such degenerate hexagons (another way of putting this is that the valence $4$, $5$, and $6$ circle patterns, which are the ones with degenerate or overlapping hexagons,
are rigid) and the different parabolic cycles can be verified by hand for all of these. As we will see in the $ \delta = 2\pi/3 $ example, it is possible for order $2$
elliptic elements to appear when there are complete overlaps: if $ C_1 $ is equal to a translate $ \lambda + C_2 $ (where $\lambda \in \Lambda$)
then the element which maps $ C_1 $ to $ \lambda + C_2 $ must be an involution.

The following example will exhibit discrete groups which come from circle patterns with all four possible tangency graph valencies.
\begin{ex}[Lattice circle patterns]\label{ex:circles}
  We have just seen that, given circle patterns in the plane satisfying the conditions (CP1), (CP2), and (CP3), we obtain a single family of discrete groups of
  the form $ G(\alpha,\beta,1) $ with fundamental domain arising from the exterior of the discs in the patterns; these groups are parameterised by the dihedral
  angle $ \delta $ of the hexagon with vertices at $ \pm 1 $, $ -1 + \alpha $, $ 1-\beta $, $ 1-\beta+\alpha $, and $ -1+\alpha-\beta $. The parameter space of
  this family is not a connected set: for instance when $ \delta $ increases slightly past $ \pi/2 $ condition (CP2) is violated and is only satisfied again
  once the angle has increased sufficiently for all overlapping circles to coincide.

  \begin{figure}
    \begin{subfigure}[t]{0.49\textwidth}
      \centering
      \includegraphics[width=\textwidth]{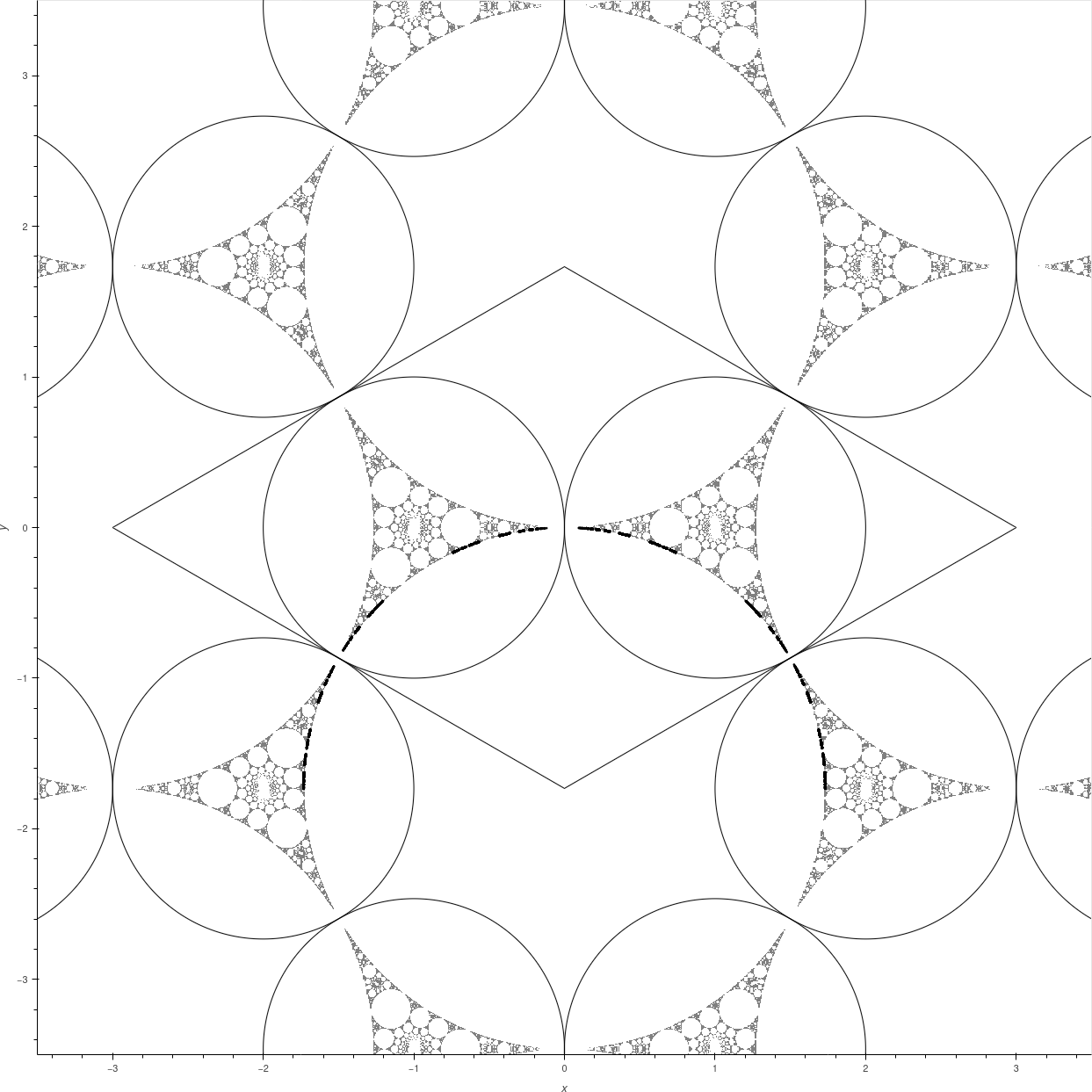}
      \caption{$ \delta = \pi/3 $\label{fig:circles_pi3}}
    \end{subfigure}\hfill%
    \begin{subfigure}[t]{0.49\textwidth}
      \centering
      \includegraphics[width=\textwidth]{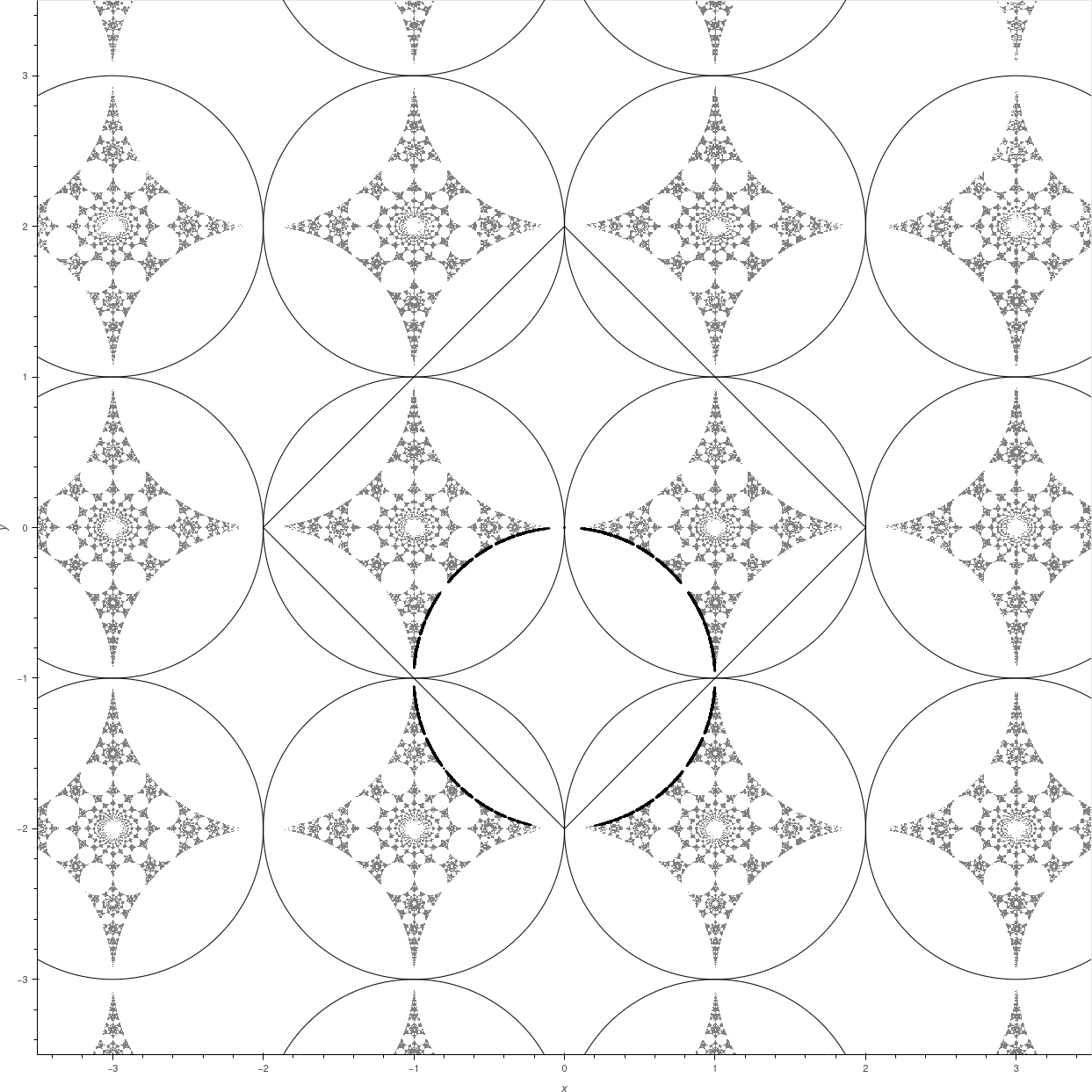}
      \caption{$ \delta = \pi/2 $\label{fig:circles_pi2}}
    \end{subfigure}\\[2ex]
    \begin{subfigure}[t]{0.49\textwidth}
      \centering
      \includegraphics[width=\textwidth]{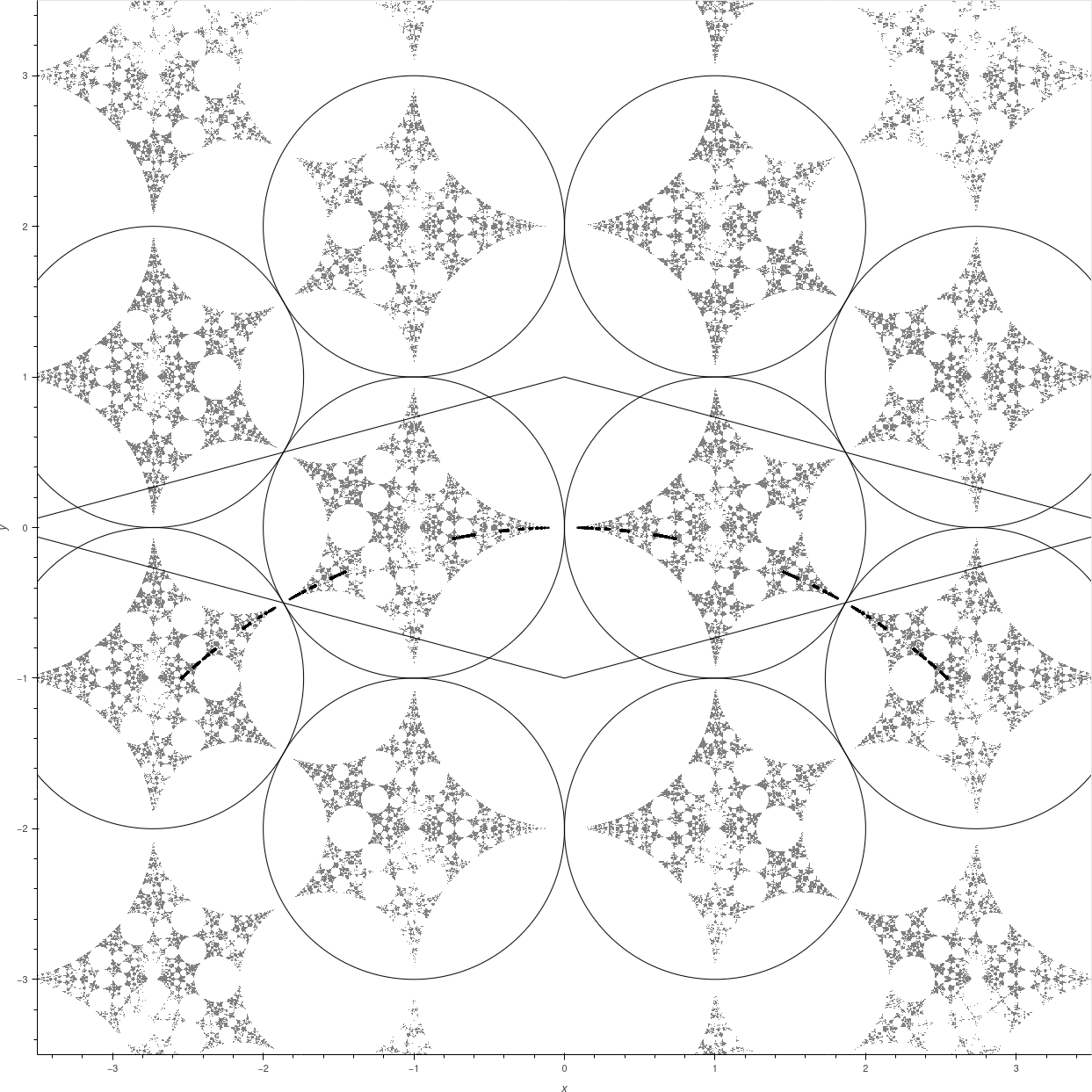}
      \caption{$ \delta = \pi/6 $\label{fig:circles_pi6}}
    \end{subfigure}\hfill%
    \begin{subfigure}[t]{0.49\textwidth}
      \centering
      \includegraphics[width=\textwidth]{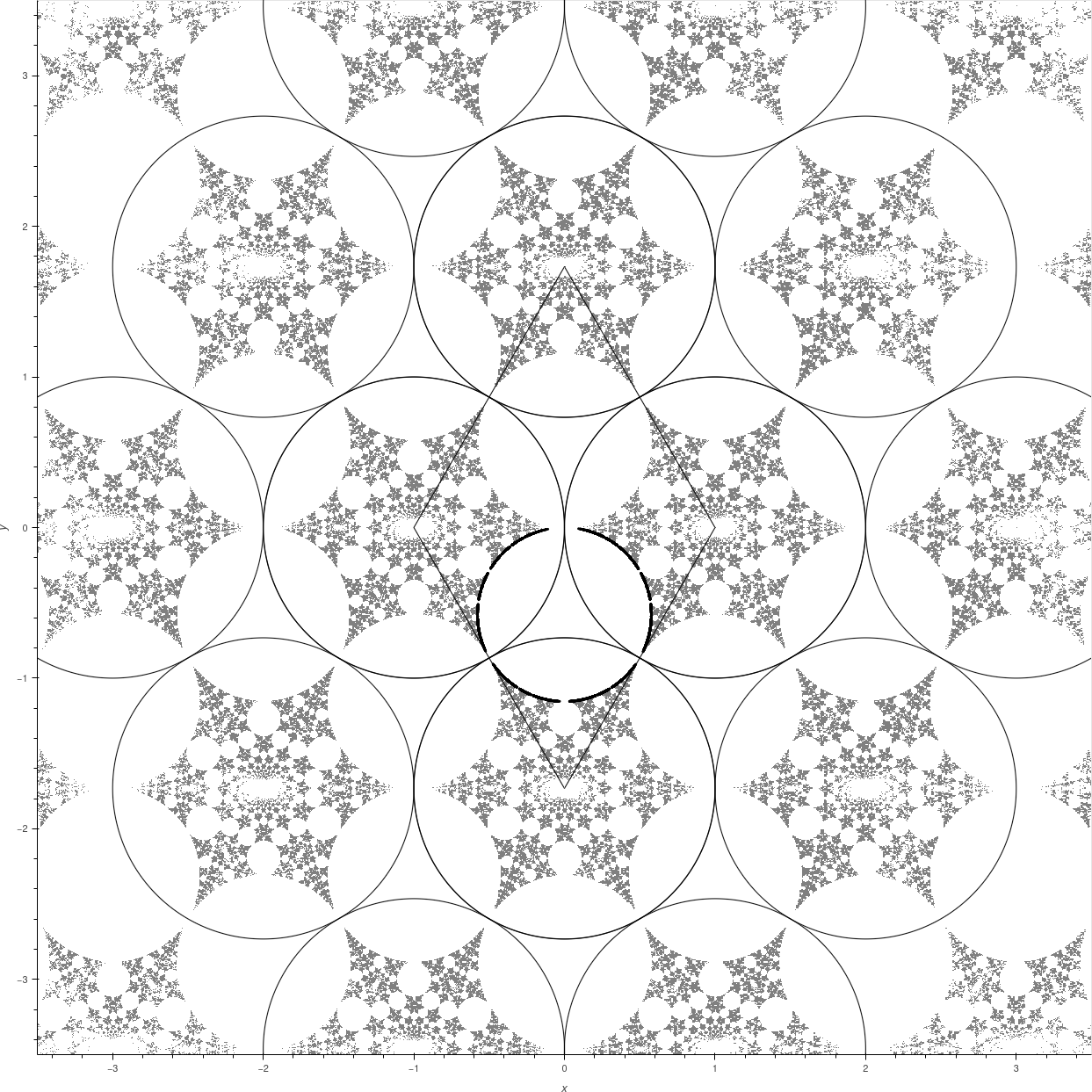}
      \caption{$ \delta = 2\pi/3 $\label{fig:circles_2pi3}}
    \end{subfigure}
    \caption{Limit sets of groups realising all four circle pattern valencies. The black points are the limit points of $ P_3 $.}
  \end{figure}

  The knowledge we have of the parabolic cycles suggests that an important subgroup of this family of groups is the subgroup $P_3$ generated
  by $ M^{-1} $ and $ MP^{-1} MQ^{-1} $. The product $ P^{-1} M Q^{-1} $ has trace $ -2-4 \cos\delta $, which is real for all $ \delta $, so $ P_3 $ is always Fuchsian.
  \begin{itemize}
    \item When $ \delta = \pi/3 $, \cref{fig:circles_pi3}, we have a regular hexagonal lattice and each circle is tangent to $3$ others. The
          element $ P^{-1} M Q^{-1} $ has trace $ -4 $, so is hyperbolic. The group $ P_3 $ is contained within a larger Fuchsian group that
          stabilises an open disc in the domain of discontinuity of $ G(\alpha,\beta,\lambda) $.
    \item When $ \delta = \pi/2 $, \cref{fig:circles_pi2}, we have a regular square lattice and each circle is tangent to $4$ others. The element $ P^{-1} M Q^{-1} $
          has trace $ -2 $, and is the vertex cycle element corresponding to the new tangency between ``opposite" circles in the hexagon. The group $ P_3 $ is
          an $ (\infty,\infty,\infty)$-triangle group which stabilises an open disc in the domain of discontinuity of $ G(\alpha,\beta,\lambda) $.
    \item When $ \delta = \pi/6 $, \cref{fig:circles_pi6}, each circle is tangent to $5$ others. One of the two additional parabolic elements that come from the vertex
          cycles at the tangencies of the ``top and bottom'' circles in the hexagon is $ PQ^{-1}MP^{-1}QM^{-1} $. The group $ P_3 $ is Fuchsian but, unlike
          the other examples here, does not preserve an open disc in the domain of discontinuity.
    \item When $ \delta = 2\pi/3 $, \cref{fig:circles_2pi3}, we have a regular triangular lattice and each circle is tangent to $6$ others. The
          element $ P^{-1} M Q^{-1} $ has trace $ 0 $, so is a rotation of order $2$. It has fixed points at $ \pm i - i\sqrt{3} $ (these are the two
          intersection points with the imaginary axis of one of the grey circles in the figure, and the element acts to exchange
          the interior and exterior of this circle). This group will become more important in \cref{ex:whitehead}.
  \end{itemize}
\end{ex}

The hexagonal structure of \cref{fig:circles_2pi3} reminded us of a series of diagrams in \autocite[\S III.5]{kag}. This motivated us to consider
certain examples of the groups $ G(\alpha,\beta,\lambda) $ which arise from representations of links.

\begin{ex}[The Whitehead link]\label{ex:whitehead}
  The group $ G $ generated by
  \begin{displaymath}
    P=\begin{bmatrix} 1 & 2 \\ 0 & 1 \end{bmatrix}\!,\;
    Q=\begin{bmatrix} 1 & i \\ 0 & 1 \end{bmatrix}\!,\;\text{and}\,
    M=\begin{bmatrix} 1 & 0 \\ -1-i & 1 \end{bmatrix}
  \end{displaymath}
  is discrete and $ \H^3/G $ is homeomorphic to the complement of the Whitehead link \autocite{wielenberg78} (also \autocite[Example~59]{kag}). The group admits a presentation
  \begin{displaymath}
    \langle P,Q,M : M P Q^{-1} M Q^{-1} M^{-1} Q M^{-1} Q^{-1} = [P,Q] = 1 \rangle.
  \end{displaymath}
  Since $ [P,Q] = 1 $, this group has a representative in the space of groups $ G(\alpha,\beta,\gamma) $ which we have been studying. Conjugating by the map $ z \mapsto (-1-i)z $,
  we find that this representative is $ G_{\mathrm{Wh.}} = G(-2-2i, 1-i, 1) $. In addition to the three generators, the element $ Q^{-1} MP Q^{-1} $ is parabolic; it represents
  a meridian of the link and is basepoint-free isotopic to the meridian represented by $ M $.

  \begin{figure}
    \centering
    \includegraphics[width=3.5in]{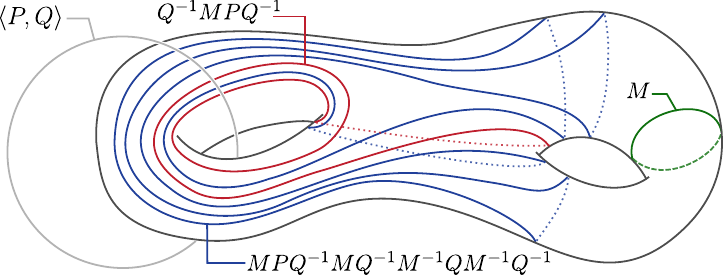}
    \caption{The curve arising from the Whitehead relator and two other curves forming a pair of pants decomposition of $ S_2 $.\label{fig:whitehead}}
  \end{figure}

  If $ \Sigma $ is the genus two surface and if $ (\alpha,\beta,\lambda) \in \C^3 $, then there is a natural map $ \pi_1(\Sigma) \to G(\alpha,\beta,\lambda) $ coming from inclusion of the
  surface into the boundary of a $ (1;2)$-compression body and then projecting down to a specific representation of the compression body holonomy group.
  It so happens that the long relator in the Whitehead group presentation lies in this image, as do the curves represented by $ M $ and $ Q^{-1} MP Q^{-1} $: see \cref{fig:whitehead}, which shows
  that the three curves give a pair of pants decomposition of $\Sigma$. These three words have trace $ \pm 2 $ in $ G_{\mathrm{Wh.}} $, and so the parameters for $ G_{\mathrm{Wh.}} $ satisfy
  the system of equations
  \begin{displaymath}
    \tr^2 M P Q^{-1} M Q^{-1} M^{-1} Q M^{-1} Q^{-1} = \tr^2 M = \tr^2 Q^{-1} MP Q^{-1} = 4.
  \end{displaymath}
  In general we refer to the solutions of a system of equations constructed in this way as the \df{parabolic locus} of the triple of curves. Despite the name it is possible
  for words to be killed rather than made parabolic, as in $G_{\mathrm{Wh.}}$.

  Solving the system to find the other representations in the parabolic locus of the three curves, we find that exactly four solutions give groups $ G(\alpha,\beta,\lambda) $ which are discrete, non-Fuchsian,
  and non-elementary. The four solutions come in pairs which are symmetric and correspond to different generating sets of the same group. One of the two groups obtained is the Whitehead link group $ G_{\mathrm{Wh.}} $.
  The other is the group $ G(1+i\sqrt{3}, -2+2i\sqrt{3}, 1) $, which is index two in the group of \cref{ex:circles} with $ \delta = 2\pi/3 $; thus it has the same circle-packing limit set. The quotient
  surface of this group is a pair of thrice-punctured spheres.
\end{ex}

\begin{figure}
  \centering
  \includegraphics[width=3.5in]{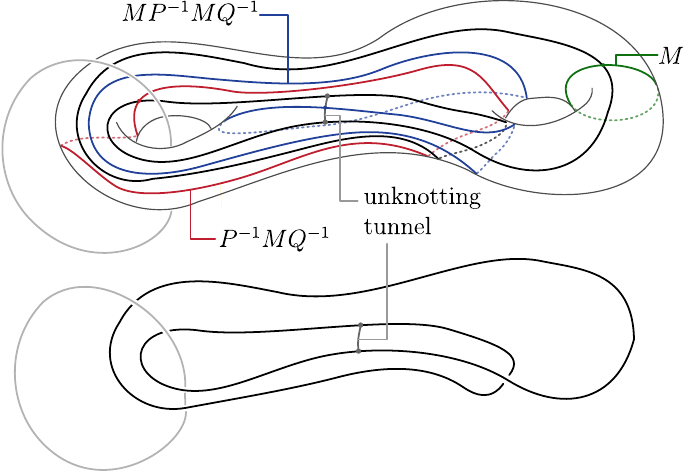}
  \caption{Solomon's knot (\href{https://katlas.org/wiki/L4a1}{\texttt{L4a1}}) is the link obtained by filling in a rank $1$ cusp in the $ \delta = \pi/2 $ group of \cref{ex:circles}. It is non-hyperbolic since
  it embeds on a torus. We will provide a different, representation-theoretic, proof of this fact later on (\cref{ex:solomon_slice}).\label{fig:solomon}}
\end{figure}

\begin{ex}[Solomon's knot]\label{ex:solomon}
  We are led to ask whether the torsion-free group $G$ with $ \delta = \pi/2 $ of \cref{ex:circles} arises from a knot or link group representation in the same way. In this
  group, $ P^{-1} M Q^{-1} $, $ M $, and $ M P^{-1} M Q^{-1} $ are parabolic. If we perform a filling on $ \H^3/G $ by gluing in a disc with boundary represented by $ M P^{-1} M Q^{-1} $,
  then the resulting manifold is a link complement. The link itself can be found by carefully drawing in the dual $\theta$-curve to the three curves on the genus $2$ surface,
  pushing it slightly into the handlebody $H$ that the surface bounds, deleting the unknotting tunnel (dual to $ MP^{-1} MQ^{-1} $), and then deleting the surface leaving only the
  rank $2$ cusp $ \langle P, Q \rangle $ and the new curve in $ \R^3 $: see \cref{fig:solomon}. This procedure arises from the theory in Cho and McCullough \autocite{cho07}; in that paper, it is shown
  that the usual definition of unknotting tunnel is equivalent to one involving embedded discs in $H$, and we are simply applying this new definition while also keeping careful track
  of algebraic and combinatorial information about the boundaries of the discs. \textit{A priori} this only works when the dual curves to the knot or link bound embedded discs
  in the handlebody $H$ since in order to get a link complement in the sphere we need to be able to push the surface curves onto a point in the bounded handlebody, and so the careful
  reader may be slightly worried to observe that the curves we have drawn in \cref{fig:whitehead,fig:solomon} do not bound discs in this way. However, because we have drilled out the
  rank $2$ cusp we are allowed to work modulo Dehn twists on the surface around loops which are homotopic to curves on the cusp torus. This is because these loops bound embedded punctured
  discs and so the Dehn twists may be undone via isotopy in the ambient manifold (c.f.\ McCollough and Miller \autocite[vii]{mccolloughmiller}). One can check that performing a Dehn twist around the meridian
  of the torus end as shown in both pictures does not change the curve that represents the link component on the surface, but does untwist the dual curves enough that they bound discs in $H$.
\end{ex}

Given any maximal system of non-intersecting closed curves on the genus two surface that are homotopically nontrivial and not homotopic to the rank two cusp---i.e.
a maximal rational lamination---we can read off three words and therefore three trace equations that we can solve to produce intricate Kleinian groups. These groups
do not, in general, come from representations of groups defined in terms of circle patterns meeting conditions (CP1)--(CP3); for instance, $ M $ is in general
non-parabolic and so $ \lambda \neq 1 $. We give an illustrative example now.

\begin{ex}[The $8_5$ knot with an extra loop]\label{ex:eightfive}
  We will reverse the procedure from the previous examples, and start from a knot in order to construct systems of three disjoint curves on the genus $2$ surface
  and hence interesting representations of $ (1;2)$-compression body groups. To do this we take a tunnel number $1$ knot, embed it on the genus $2$ end of the topological
  compression body, and use the combinatorics of this embedding to define words in the group onto which we impose trace conditions. If this is done with care then we obtain groups related
  to the holonomy group of the tunnel number $1$ knot together with an additional unknotted loop. In \cref{ex:eightfive_slice} below we will perform the same procedure for the knot directly in a
  genus $2$ handlebody, without the additional loop.

  \begin{figure}
    \centering
    \begin{subfigure}[m]{0.35\textwidth}
      \centering
      \includegraphics[width=\textwidth]{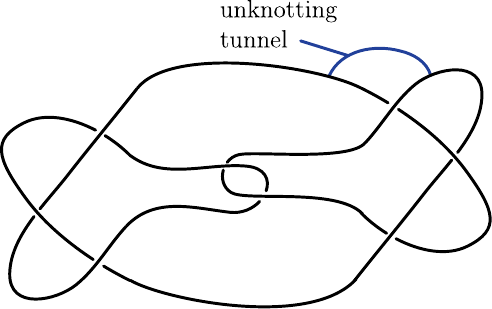}
      \caption{A knot diagram for $ \mf{p} $, showing an unknotting tunnel.\label{fig:eightfive_knot}}
    \end{subfigure}\hfill%
    \begin{minipage}[m]{0.6\textwidth}
    \begin{subfigure}{\textwidth}
      \centering
      \includegraphics[width=\textwidth]{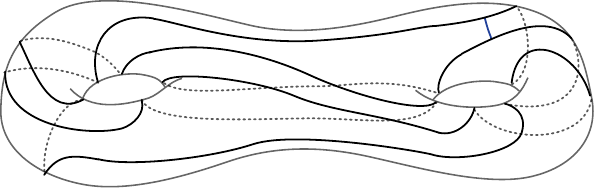}
      \caption{Embedding of $ \mf{p} $ on the genus $2$ surface.\label{fig:eightfive_surface}}
    \end{subfigure}\\[2ex]
    \begin{subfigure}{\textwidth}
      \centering
      \includegraphics[width=\textwidth]{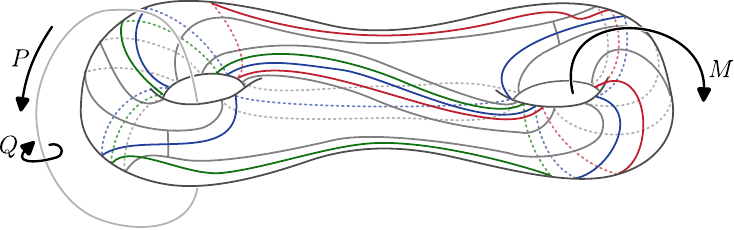}
      \caption{Two meridian curves (red and green) and one curve dual to the unknotting tunnel (blue). The labelled arrows indicate the direction that the curve coding for the $ (1;2)$-compression body
      representation was done (e.g.\ if a curve winds backwards along the arrow labelled $ P $, then its word picks up a letter $P^{-1} $).\label{fig:eightfive_dual}}
    \end{subfigure}
    \end{minipage}
    \caption{The $8_5$ knot $ \mf{p} $.}
    \end{figure}

  One hyperbolic knot of minimal crossing number which has tunnel number $1$ but which is not $2$-bridge is the $ 8_5 $ knot $ \mf{p} $, also known as the $ (3,3,2) $ pretzel knot. The symmetric
  presentation of $ \mf{p} $ given in \cref{fig:eightfive_knot} exhibits an unknotting tunnel $ \tau $ and also makes obvious an embedding onto the genus $2$
  surface, \cref{fig:eightfive_surface}. Since the knot has tunnel number $1$ there exists an embedding in which the dual curves to the $\theta$-curve $ \mf{p} \union \tau $ bound discs, but
  for this knot the words involved would get very long (more than $30$ letters) and the diagrams of the curves on the surfaces become much more intricate and accordingly
  harder to read. Our embedding is simpler to work with but will still produce interesting groups, and the reader interested in the knot theory can write the correct words
  down and follow the same steps to produce the link complement. An explicit algorithm, and a worked example for the twisted torus knot $ T(5,7,2) $, can be found in Ishihara \autocite{ishihara11}.

  In \cref{fig:eightfive_dual}, we show the curves which are dual to the three arcs of the $\theta$-curve $ \mf{p} \union \tau $ with respect to our embedding. Drilling out a rank two cusp
  from the handlebody to view the surface as the genus $2$ surface in a $(1;2)$-compression body and reading out the words of the dual curves in terms of the generating
  set $ \{P,Q,M\} $ gives us the three words listed in \cref{tab:words}.

  \begin{table}
    \caption{Words arising dual to the embedded $ 8_5 $ knot.\label{tab:words}}
    \begin{tabular}{ll}\toprule
      $ U_1 = P^{-1}MM $ & (red)\\
      $ U_2 = MMMP^{-1}P^{-1}QP^{-1} $ & (unknotting tunnel dual, blue)\\
      $ U_3 = MP^{-1}P^{-1}Q $ & (green)\\\bottomrule
    \end{tabular}
  \end{table}

  There are two solutions (modulo various symmetries) to the system $ \tr^2 U_1 = \tr^2 U_2 = \tr^2 U_3 = 4 $ with $ \alpha $ and $ \beta $ not linearly dependent over $ \R $, namely
  \begin{displaymath}
    (\alpha,\beta,\lambda) =\begin{cases}
      \left( \frac{1}{2}(3-i\sqrt{7}), \frac{1}{2}(7-i\sqrt{7}), \frac{1}{4}(3-i\sqrt{7}) \right)\;\text{or}\\
      (1.7581 - 2.7734i, 6.4537 - 4.8311i, -0.4688 - 0.3578i).
    \end{cases}
  \end{displaymath}

  \begin{figure}
    \begin{subfigure}[t]{0.49\textwidth}
      \centering
      \includegraphics[width=\textwidth]{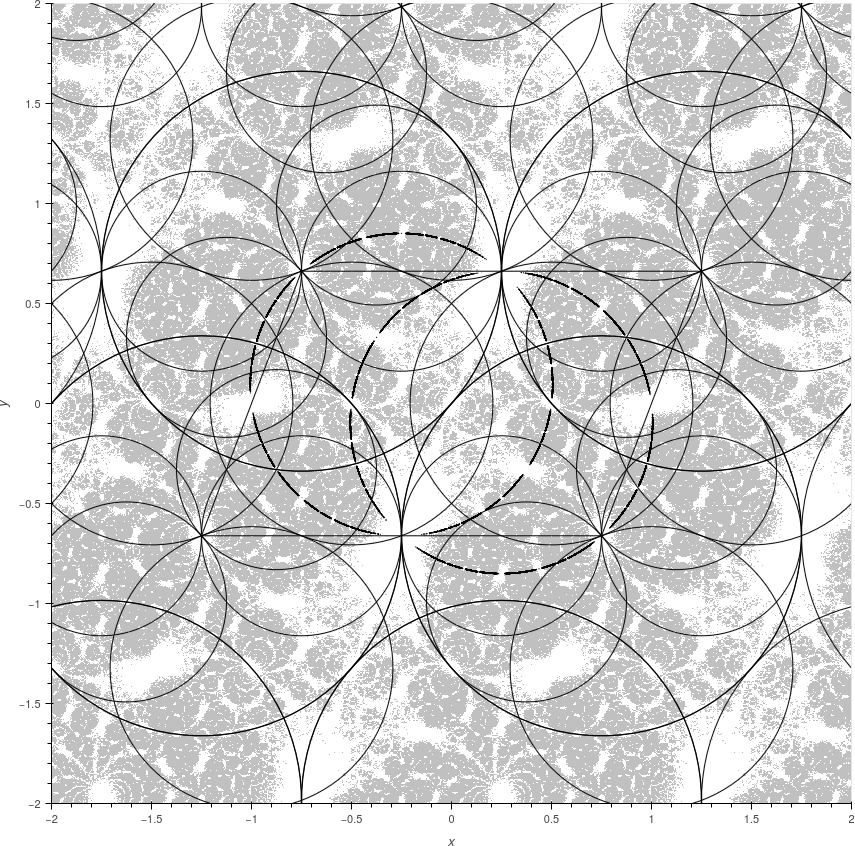}
      \caption{The cofinite volume group. Isometric circles of elements related to the combinatorics of the embedding of the knot are shown.\label{fig:eightfive_lim_knot}}
    \end{subfigure}\hfill%
    \begin{subfigure}[t]{0.49\textwidth}
      \centering
      \includegraphics[width=\textwidth]{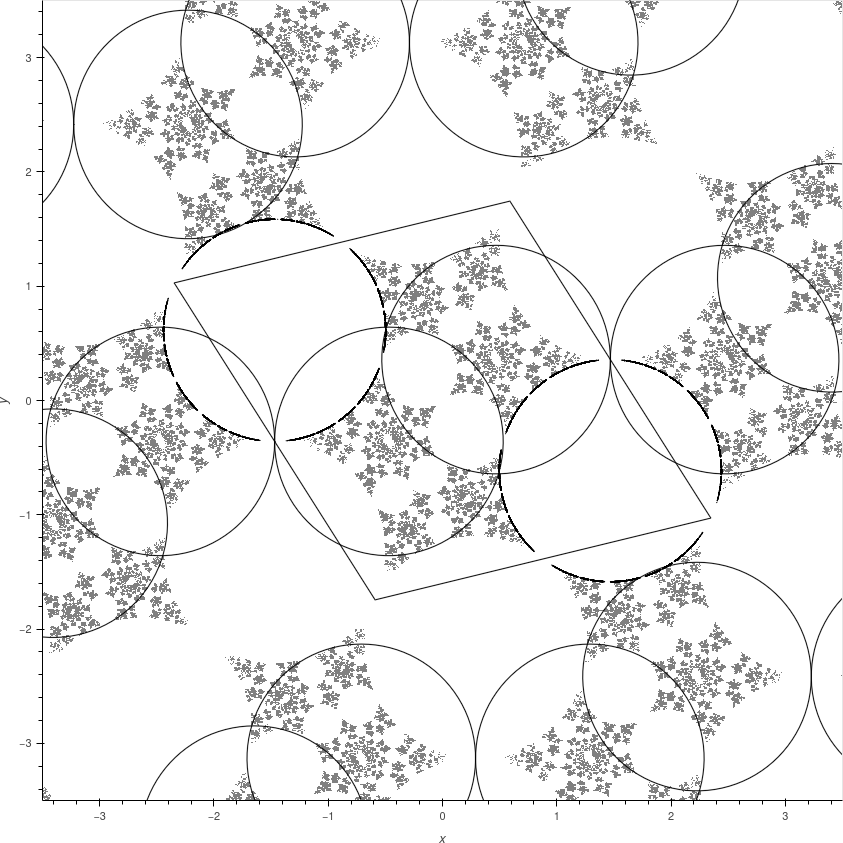}
      \caption{The group on the boundary of compression body space with circle-packing limit set. This example does not arise from a circle pattern like those above: the isometric
      circles of $ M $ (shown) overlap each other.\label{fig:eightfive_lim_cusp}}
    \end{subfigure}
    \caption{Groups arising from a `bad' genus $2$ embedding of the $ 8_5 $ knot $ \mf{p} $, with an additional torus drilled out.}
  \end{figure}

  The first listed solution gives a group $ G_1 $ with $ \Lambda(G_1) = \hat{\C} $; the quotient $ \H^3/G_1 $ is finite volume. We show a small number of points from the limit set of $ G_1 $,
  overlaid with some combinatorial data, in \cref{fig:eightfive_lim_knot}. Since our curves did not bound discs, $G_1 $ is \emph{not} the holonomy group of the complement of the $ 8_5 $ knot
  with an additional loop drilled out; the problem that the curves do not bound discs can't be fixed by Dehn twists as was possible in \cref{ex:solomon}, because the twists
  that would be needed modify the knot curve as well as the dual curves. The group is not even torsion-free, since $ \Phi = M^{-1}P $ is
  an elliptic of order $2$. The word $ U_1 $ is conjugated by $ \Phi $ to $ P $ and hence the loop which it represents is isotopic to the meridian of the cusp uniformised by $ \langle P, Q \rangle $;
  there are two rank $2$ cusps, $ \langle P,Q \rangle \cong \langle U_1, MP^{-1}QMP^{-1} \rangle $ and $ \langle M^{-1} U_2 M, U_3 \rangle $, and one closed elliptic cone arc with longitude and meridian lying
  in the elementary subgroup $ \langle MP^{-1}, Q^{-1}P^2M^{-2}Q^{-1}P^3 \rangle $.

  A formal proof of discreteness of $ G_1 $ can be obtained from the Poincar\'e polyhedron theorem; we neglect the details, but overlaid over the limit set in the figure we have drawn
  a parallelogram fundamental domain for the lattice $ \langle P, Q \rangle $ (the edges of the parallelogram are paired by $ QP^{-1} $ and $ P^2 Q^{-1} $) together
  with the isometric circles of $ U_1 $, $ U_2 $, and $ \Phi Q \Phi $, and translates of these circles under the lattice. The domes over these circles can be used to construct a fundamental
  polyhedron in $ \H^3 $ for the action of the group.

  The second solution to the system of equations, which is cubic over $ \Q $, gives a group $ G_2 $ with no rank $2$ parabolic subgroups other than $ \langle P,Q \rangle $; we will later
  refer to this group as the cusp group arising from the system of equations. The conformal boundary $ \partial \H^3/G_2 $ consists of a pair of thrice-punctured spheres, joined by
  the three rank $1$ parabolic cusps represented by $ \langle U_1 \rangle $, $ \langle U_2 \rangle $, and $ \langle U_3 \rangle $. The manifold $ \H^3/G_2 $ can be viewed as the limit
  of a sequence of submanifolds of $ \Sph^3 $ obtained by taking one handlebody of a genus $2$ Heegaard splitting and allowing three curves on the Heegaard surface to continuously contract
  in length by shrinking into the second `deleted' handlebody. Since the curves do not bound discs in the second handlebody, this sequence of embeddings does not limit onto an embedding
  of $ \H^3/G_2 $ into $ \Sph^3 $ even though the limit manifold is a compression body manifold.

  We show $ \Lambda(G_2) $ in \cref{fig:eightfive_lim_cusp}: the complement of the limit set is a circle packing,
  and the circles shown are the isometric circles of $ M $ (together with lattice translates). Additionally,
  we show in black the limit sets of a pair of embedded Fuchsian thrice-holed sphere subgroups $ \langle U_1,U_3 \rangle $ and $ \langle U_1^\dagger, U_3^\dagger \rangle $ (here, $ (-)^\dagger $ denotes
  reversal of the word). These two subgroups are non-conjugate and uniformise the two hyperbolic thrice-punctured spheres in the boundary of the convex core of the quotient $3$-manifold $ \H^3/G_2 $. The same two
  pairs of words also define Fuchsian subgroups of $ G_1 $ that uniformise intersecting thrice-punctured sphere subsurfaces of $ \H^3/G_1 $.
\end{ex}

In the previous example some insight was needed to find the correct group elements to produce a fundamental domain. For a discussion of automatic calculation of Ford domains for knot groups via a `greedy'
algorithm, see Riley \autocite{riley82,riley83}. More efficient automatic computations are possible if additional structural information is known: for instance, Page \autocite{page15} has described an
algorithm to find fundamental domains of arithmetic groups. Fundamental domains for infinite covolume groups are also well-studied and the most well-developed algorithms are for punctured torus groups,
implemented by Wada \autocite{wada06} and described by Yamashita \autocite{yamashita12}.

In this section we have seen several examples of groups with circle-packing limit sets. There exists a generous supply of these groups since every maximal curve system on a genus $2$ surface can
be used to obtain a circle packing with a $(1;2)$-compression body symmetry group; this follows from general results of Keen, Maskit, and Series~\autocite{keen91} for function groups.
Each of these circle packing groups gives rise to a lattice $ \alpha \Z \oplus \beta \Z $, which is called the \df{shape of the cusp} of the corresponding 3-manifold. Dang and Purcell~\autocite{dang19}
posed the following inverse problem:
\begin{problem}\label{prob:inverse_problem}
  Classify all lattices $ \alpha \Z \oplus \beta \Z $ which can arise from a maximal system of rank one parabolics on the genus $2$ end of a $ (1;2)$-compression body. For each such
  lattice give a construction of the corresponding group.
\end{problem}
This is part of a family of problems dating back at least to Nimershiem \autocite{nimershiem94} who proved that, arbitrarily close to any complex structure on the torus, there exists
a complex structure admitting an embedded circle pattern that lifts to a pattern of tangent circles on the Riemann sphere and so appears as a lattice in some hyperbolic 3-manifold
group: in other words, a dense selection of lattices can appear as cusp shapes. Nimershiem's theorem relies on general machinery of Brooks \autocite{brooks86} who used the lifting of circle patterns from
analytically finite Riemann surfaces to the Riemann sphere to prove that Kleinian groups admitting cofinite extensions are dense in deformation spaces of geometrically finite groups.

An elementary introduction to the geometry of \cref{prob:inverse_problem} in a similar vein to our exposition here may be found in \autocite[Chapter~9]{indras_pearls}, and Wright~\autocite{wright05}
extended this work to study some versions of this inverse problem for cusp groups on the boundary of the Maskit and Riley slices and for maximal cusps of higher genus Schottky spaces. The modern study
of \cref{prob:inverse_problem} involves some very deep number theory, and a very accessible introduction is Stange \autocite{stange24}. Of additional interest in this direction is recent work of Kapovich and
Kontorovich~\autocite{kapovich23} on groups obtained by taking Kleinian groups with circle-packing limit sets, or more generally limit sets with patterns of intersecting discs, and extending them
by adjoining the reflections in all the discs in the circle packing.

\section{Linear sections of character varieties}\label{sec:families}
It would be nice to be able to visualise entire spaces of representations. In the literature such visualisations do exist, for example for the \df{Maskit slice} \autocite[Fig.~9.11]{indras_pearls}
and Bers slices \autocite{komori06} of once-punctured torus groups, and for the \df{Riley slice} of four-punctured sphere groups (see \cref{defn:rileyslice} below) \autocite[Fig.~1]{ems21}. We
will present here a general method for drawing such pictures that does not rely on the special trace relations or other combinatorial structures that exist in those settings. Our method is very
na\"ive, and is of interest only because it is easy to implement, is surprisingly fast in practice, and does not require any specialised analysis for applying it in different deformation spaces.

We will begin with an illustrative example of a character variety parameterisation and corresponding slice.
\begin{defn}
  A \df{Schottky group} of genus $ g $ is a purely loxodromic discrete subgroup of $ \PSL(2,\C)$ that is free on $ g $ generators; the corresponding hyperbolic
  manifold is a handlebody of genus $g$. The locus of such groups is a subset of the character variety $ \Hom(F_g, \PSL(2,\C)) \git \PSL(2,\C) $ called
  \df{Schottky space}, $ \mc{S}_g $. It is a very complicated open set, and admits a complex structure of dimension $ 3g - 3 $ as a quotient of the Teichm\"uller
  space of genus $g$ surfaces \autocite[Theorem~5.27]{matsuzakitaniguchi}. A \df{maximal cusp point} of $ \mc{S}_g $ is a point on $ \partial \mc{S}_g $ which
  is still free on the induced generators from $ F_g $ and is still purely loxodromic except for $ 3g-3 $ rank one parabolic subgroups, the maximal possible number of such subgroups.
\end{defn}

We will not explore this direction in the current paper, but Schottky space has strong connections to the computational geometry of Riemann surfaces and to the moduli of algebraic
curves; see for example Chapter~3 of Bogatyrev \autocite{bogatyrev} and recent work of Fairchild and Ríos Ortiz \autocite{fairchild24}. This point of view has found applications
in mathematical physics \autocite{ichikawa25}.

\begin{defn}\label{defn:grandma}
  We introduce a parameterisation of $ \mc{S}_2 $, following Horowitz \autocite{horowitz72} but adapting his normalisation to include some of the features of
  the \df{Grandma's recipe} normalisation of punctured torus space given in \autocite[Chapter~8]{indras_pearls}. Set
  \begin{displaymath}\renewcommand*{\arraystretch}{1.5}
    X = \begin{bmatrix} \frac{1}{2}(t_X + i t_{XY}) & -\frac{1}{2} (t_X + v) \\ -\frac{1}{2} (t_X - v) & \frac{1}{2}(t_X - i t_{XY}) \end{bmatrix}\;\text{and}\;
    Y = \begin{bmatrix} \frac{t_Y}{2} - i & \frac{t_Y}{2} \\ \frac{t_Y}{2} & \frac{t_Y}{2} + i \end{bmatrix}
  \end{displaymath}
  where $ v $ satisfies $ v^2 = 4 - t_{XY}^2 $. Then $ \tr X = t_X $, $ \tr Y = t_Y $, and $ \tr XY = t_{XY} $. By abuse of notation, we write $ \mc{S}_2 $ for the
  subset of $ \C^3 $ consisting of $ (t_X, t_Y, t_{XY}) $ for which the group $ \langle X, Y \rangle $ is Schottky of genus $2$.
\end{defn}

We will need the function
\begin{displaymath}
  \tau(M) = \abs{2\log \frac{\abs{\tr M + \sqrt{(\tr M)^2-4}}}{2}}
\end{displaymath}
which sends $ M \in \PSL(2,\C) $ to its real translation length \autocite[Exercise 7-20]{marden}; $ \tau(M) = 0 $ if $ M $ is elliptic, parabolic, or the identity.

\begin{figure}
  \centering
  \includegraphics[width=4in]{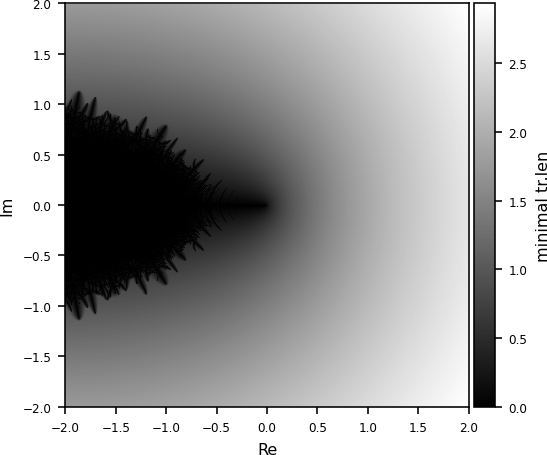}
  \caption{A linear slice through $ X(F_2) $ with $ \tr X = \tr Y = \tr XY $. The light region is an approximation to the intersection
          of the slice with $ \mc{S}_2 $, and the dark region contains indiscrete groups and sporadic discrete groups (e.g.\ knot complements and related orbifolds).\label{fig:schottky_slice_easy}}
\end{figure}

\begin{ex}[{The diagonal slice of $ \mc{S}_2 $ \autocite{series17}}]\label{ex:diagonal}
  With the normalisation in \cref{defn:grandma}, the point $(-2,-2, -2) \in \C^3 $ lies on the boundary of $ \mc{S}_2 $ and $ z_0 = (-3,-3,-3) $ lies
  in the interior. These inclusions can both be proved by drawing Ford domains. In the matrices of \cref{defn:grandma} we substitute $ -3t + -2(1-t) = 5t-2 $ for all three of $ t_X $, $ t_Y $,
  and $ t_{XY} $, giving a linear slice of $ X(F_2) $ parameterised by $ t \in \C $ such that $ t = 0 $ corresponds to the cusp point and $ t = 1 $ corresponds to $ z_0 $.
  The parameterisation induces a substitution map $ \rho : \C \to \PSL(2,\C)^2 $, which sends $ t $ to the pair of matrices $ (X,Y) $ with the parameters
  substituted---formally, $\rho$ is induced by the composition of $ t \mapsto (5t-2,5t-2,5t-2) $ with the canonical rational map $ \C^3 \times \PSL(2,\C(t_X,t_Y,t_{XY})) \to \PSL(2,\C) $.
  Within each group $ \langle \rho(t) \rangle $, we can compute the minimal value of $ \tau $ over all words up to some fixed length. For values of $ t $ such that this
  group lies in $ \mc{S}_2 $, the minimal value of $ \tau $ will be strictly positive. Conversely, if there is a nontrivial word in $ F_2 $ descending to a word $ w \in \langle \rho(t) \rangle $
  with $ \tau(w) = 0 $ then the corresponding group $ \langle \rho(t) \rangle $ cannot lie in Schottky space. By plotting the minimal value of $ \tau $ in $ \langle \rho(t) \rangle $ for each $ t $, we
  obtain \cref{fig:schottky_slice_easy}. The big white/light grey region corresponds to discrete groups with isomorphism type $ F_2 $ that uniformise handlebodies,
  i.e.\ it is a slice through $ \mc{S}_2 $ (or rather, an approximation to it: the boundary, which is highly fractal, is not seen to a very high resolution).
\end{ex}

The general algorithm, which works for arbitrary groups and not just free ones, is not much more complicated. Following the algorithm, we will say a few words about the qualitative
and quantitative meanings of the image which it produces.

\begin{alg}[Drawing character variety slices]\label{alg:stupid_picture}
  Let $ \rho : \C \to \PSL(2,\C)^g $ be a holomorphic map, which should be viewed as defining a 1-parameter family of subgroups of $ \PSL(2,\C) $ by giving
  their generators, just as in \cref{ex:diagonal} above. Often, as in that example, we will define $ \rho $ by taking a fixed set of matrices depending on some
  parameters (in the example, $ \{X,Y\} $ which depended on $ t_X$, $t_Y$, and $t_{XY} $) and substituting for those parameters some fixed function of a single
  parameter $ t $ (in the example, this was the function $ t \mapsto (5t-2,5t-2,5t-2 $).

  Choose a test point $ z_0 \in \C $, which gives a fixed base group $ \Gamma_0 = \langle \rho(z_0) \rangle $. This test point should be a `typical group' in the
  family: for instance to draw pictures of $ \mc{S}_2 $, $ \Gamma_0 $ should be a known Schottky group in the image of $ \rho $. For simplicity, we assume
  that $ G_0 $ is geometrically finite. Define also a set of points $ U \subset \C $ (which does not necessarily need to contain $ z_0 $) which will be the
  set of points to sample (`pixels') for the picture, and a number $ N $ which will be the maximal length of a word in the group to sample.

  After these choices, we perform the following steps:
  \begin{enumerate}
    \item Enumerate all words up to length $N$ in the group $ \Gamma_0 $, in a fixed order $ W_1,\ldots, W_n,\ldots$. One may use relations which hold in the entire
          image of $ \rho $, but not relations which apply specifically to groups in the isomorphism type of $ \Gamma_0 $, to do this, so long as the
          order is well-defined. Produce a list $ \mc{L}(z_0) $ such that the $ n$th element is $ \tau(W_n) $ for all $ n $.
    \item Now for each pixel $ z \in U $:
      \begin{enumerate}
          \item Enumerate all the words $ \{W_1,\ldots,W_n,\ldots\} \in \langle\rho(z)\rangle $ in the same order as the words in $ \Gamma_0 $, producing
                a list $ \mc{L}(z) $ of elements $ \tau(W_n) $.
          \item Compare $ \mc{L}(z_0) $ to $ \mc{L}(z) $ pairwise: if an element of the first list is $0$, then throw away the corresponding element of the second list. This
                produces a list $ \mc{L}'(z) $ of real translation lengths of elements of $\langle\rho(z)\rangle$ corresponding to words which are not elliptic, parabolic,
                or the identity at the test point $ z_0 $ and hence at every point in the domain $ U $.
          \item Assign the value $ \min \mc{L}'(z) $ to the pixel $ z $.
      \end{enumerate}
  \end{enumerate}
  The result is an array of numbers, one for each point in $ U $, that gives the minimal translation length of an element (of length $ \leq N $) of the group at $ U $,
  excluding elements which have a constant translation length of $ 0 $ everywhere in $U$.
\end{alg}
\begin{rem}[Practical improvements to the algorithm]\label{rem:improved_alg}
  The number of matrices computed explodes exponentially with the depth $ N $. However, many of the matrices produced are irrelevant (e.g.\ are conjugates of each other or are known to lie in a parabolic subgroup),
  and so if the word list is produced and trimmed using these filters the memory usage is improved. Significant improvements can also be made by choosing random subsets of words to sample instead of taking
  minima over the entire tree; one just needs to ensure that the same word list is used for the test point and all pixels.
\end{rem}

The map $ \rho : \C \to \PSL(2,\C)^g $ is a holomorphic parameterisation of a 1-dimensional subset of a character variety $ X(R) = \Hom(R,\PSL(2,\C))\git\PSL(2,\C) $, where $ R $
is the isomorphism type of the group generated by the image of $ \rho $ generically: i.e.\ $ R $ is the free group on $ g $ elements, modulo relations which
hold in $ \langle \rho(z) \rangle $ for all $ z \in \C $. More precise information on character varieties may be found in Kapovich~\autocite[\S 4.3]{kapovich}.

The fixed basepoint $ \Gamma_0 $ lies in an open subset of $ X(R) $, the \df{quasiconformal deformation space} of $ \Gamma_0 $, which is the maximal connected
open subset of $ X(R) $ containing discrete groups $ \tilde{G} $ that are isomorphic to $ G_0 $ such that (i) the isomorphism induces a homeomorphism between
the hyperbolic $3$-orbifold $\H^3/G_0$ and the hyperbolic $3$-orbifold $ \H^3/\tilde{G} $, and (ii) the two Riemann surfaces $ \Omega(G_0)/G_0 $ and $ \Omega(\tilde{G})/\tilde{G} $
lie in the same Teichm\"uller space \autocite[\S 5.3]{matsuzakitaniguchi}. The tameness theorem \autocite[\S 5.6]{marden} implies that as you move to the boundary
of the quasiconformal deformation space there must be a sequence of elements of the group which have real translation length converging to $0$, and generalisations
of McMullen's cusp density theorem \autocite[Theorem~5.3.1]{marden} tell us that there is a dense subset of points on the boundary where the minimal length is actually $ 0 $.

The set of points with length $0$ output by the algorithm definitely lie outside the deformation space or on its boundary, those with very small length may lie outside the
deformation space, and those points with very large minimal length probably lie inside the deformation space. Around the boundary the picture can get a bit
fuzzy and should only be viewed as a visual aid for what is going on. Increasing $ N $ will give increasingly better pictures, and in practice (for
simple deformation spaces that arise in the wild) we get reasonable pictures with $ N $ in the low double digits.

\begin{rem}[Comparison with the Mandelbrot set]
  The Sullivan--McMullen dictionary between the theory of Kleinian groups and the theory of iterated quadratic maps \autocite[Table~1.1]{mcmullenRFC} tells us
  that quasiconformal deformation spaces are analogous to the Mandelbrot set $M$. Pictures of $M$ can be drawn by picking a window $ U \subset \C $, a
  disc $ D $ containing $0$, and a point $ p \in \C $, and then for every pixel $ c \in U $ computing the image of $ p $ under the $k$-fold composition
  of $ z \mapsto z^2 + c $. If the image lies in $D$ then $ c $ is deemed likely to lie in $ M $. A discussion of this process together with observations
  on good choices for the number of iterations $k$ and the radius of the disc $D$ may be found in Rojas \autocite{rojas98}.

  The modified version of our algorithm described in \cref{rem:improved_alg} can be viewed completely analogously to this Mandelbrot algorithm: the role of
  the point $ p $ is taken by the identity matrix $ I_2 $, the role of the repeated composition of $ z^2 + c $ is taken by repeated multiplication of random
  elements of a generating set of the group parameterised by the point $ c $ (i.e.\ the iteration process is to at each step multiply on the left by a random
  element of $ \rho(c) $ where $ \rho $ is the function of \cref{alg:stupid_picture}), and the role of the absolute value function in determining whether
  $ c $ lies in the moduli space is instead taken up by the function $ \tau $.

  This viewpoint shows why visualising the moduli spaces of Kleinian groups is harder than visualising Mandelbrot sets. Firstly, we have additional computational
  complexity arising from working with iteration of matrix-valued functions instead of $\C$-valued functions. Secondly, instead of computing a single iterate we
  must compute many different iterates and take the minimal value of $ \tau $ over all of these choices. Finally, the number of random samplings of words required
  for each pixel in order to produce a decent picture is highly dependent on the particular parameterised generating set. If the set is chosen to be `geometrically
  simple' as in the examples throughout this section, then only tens of iterations are needed
  to produce large-scale approximations of the moduli space. On the other hand, to get fine-scale pictures of particular parts of the moduli space (like the famous
  `zooms' of the Mandelbrot set) we should choose a generating set which is adapted to detect the boundary in the desired region; that is, words which correspond
  to prominent cusp groups in that region should be of short length in terms of the chosen generators.
\end{rem}

\begin{figure}
  \centering
  \includegraphics[width=4in]{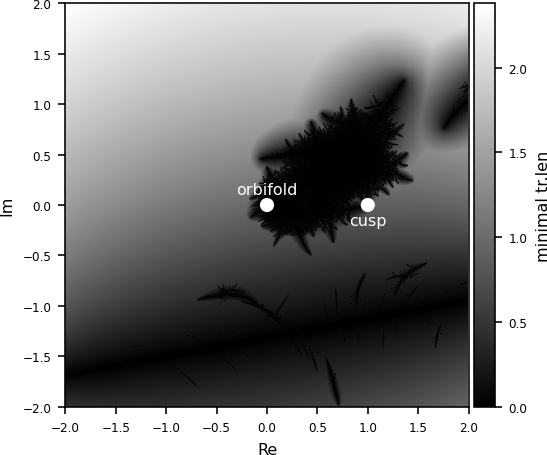}
  \caption{A linear slice through $ X(F_2) $ spanned by the cusp and orbifold group constructed from the $ 8_5 $ knot.\label{fig:schottky_slice_hard}}
\end{figure}

\begin{ex}\label{ex:eightfive_slice}
  The variety of \cref{ex:eightfive} can be mapped into the representation variety $ X(F_2) $ by killing the generator $ Q $. This involves taking the same words $ U_1 $, $ U_2 $,
  and $ U_3 $ from \cref{tab:words} but with the substitution $ P = X $, $ Q = 1 $, and $ M = Y $, and solving the system of equations
  \begin{displaymath}
    \tr^2 U_1 = \tr^2 U_2 = \tr^2 U_3  = 4
  \end{displaymath}
  for $ t_X $, $ t_Y $, and $ t_{XY} $. Up to symmetry there are two non-Fuchsian solutions:
  \begin{displaymath}
    (t_X, t_Y, t_{XY}) =\begin{cases}
      z_0 \coloneq (2-i, -i, 2-2i)\;\text{or}\\
      z_1 \coloneq (0.7607 + 0.8579i, -0.7610 - 0.8579i, 2.3146 - 2.6103i).
    \end{cases}
  \end{displaymath}
  Both solution groups are discrete: the first is a cofinite orbifold group, and the second (cubic over $\Q$) is a cusp group on the boundary of $ \mc{S}_2 $. We visualise the affine slice of the parameterisation
  $ (t_X, t_Y, t_{XY}) $ of $ X(F_2) $ that is spanned by these points: that is, we apply \cref{alg:stupid_picture} to $ t \mapsto tz_1 + (1-t) z_0 $ (as described there, this map gives a parameterisation of
  the values $ (t_X, t_Y, t_{XY}) $ which are to be substituted into the matrices of \cref{defn:grandma} to produce the family $ \rho : \C \to \PSL(2,\C)^2 $ of generating sets). The resulting picture is
  shown in \cref{fig:schottky_slice_hard}. If the embedding in \cref{ex:eightfive} is chosen differently so that the dual curves to the $ \theta$-curve $ \mf{p} \union \tau $ bound discs and the corresponding
  representations of $ F_2 $ are computed as in this example so that the dual word to $ \tau $ is killed and the other two dual words were parabolic, then one of the discrete representations obtained would
  give the unique hyperbolic structure on the $ 8_5 $ knot complement.
\end{ex}

We now move to the space of $ (1;2)$-compression bodies, which lies in $ X((\Z \oplus \Z) * \Z) $. We continue to use the normalisation $ G(\alpha,\beta,\lambda) $ of \cref{sec:lattices}.

\begin{ex}[Horizontal slice through $ (1;2)$-compression body space]\label{ex:horiz_slice}
  Consider the group $ G(\alpha,\beta, \lambda + \beta/2) $. The map $ \Phi(z) = z + \beta/2 $ acts on this group by conjugation, preserving the parabolic
  subgroup $ \langle P,Q \rangle $. We observe also that
  \begin{displaymath}
    \Phi^{-1}
    \begin{bmatrix} 1 & \beta \\ 0 & 1 \end{bmatrix}
    \begin{bmatrix} \lambda & \lambda^2 - 1 \\ 1 & \lambda \end{bmatrix}
    \Phi =  \begin{bmatrix}
                   \frac{\beta}{2}+\lambda & \left(\frac{\beta}{2}+\lambda\right)^2 - 1 \\
                   1 & \frac{\beta}{2}+\lambda
                  \end{bmatrix}
  \end{displaymath}
  which shows that the element $ M $ of $ G(\alpha,\beta, \lambda + \beta/2) $ is sent to the element $ QM $ of $ G(\alpha,\beta, \lambda) $ and
  so we see that the two groups are conjugate as groups (though the generating triples are not conjugate). Continuing the argument along these lines,
  one can show that the isomorphism type of $ G(\alpha,\beta,\lambda) $ is periodic in the $\lambda$-component by the lattice generated by $ \alpha/2 $ and $ \beta/2 $.

  \begin{figure}
    \centering
    \includegraphics[width=4in]{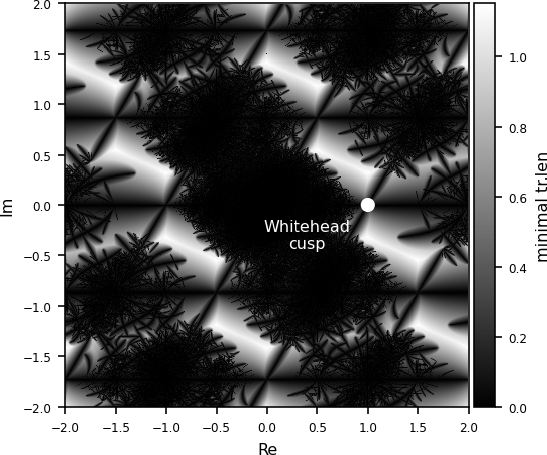}
    \caption{A `horizontal' linear slice (the $ \lambda$-plane with $ \alpha $ and $ \beta $ fixed) through $ X((\Z \oplus \Z) * \Z) $ , showing
            the periodicity that comes from the lattice.\label{fig:wielenberg_slice_horiz}}
  \end{figure}

  We can see this periodicity in \cref{fig:wielenberg_slice_horiz}. This image shows the `horizontal slice' with $ \alpha $ and $ \beta $ fixed that contains
  the cusp point $ G(2+2i\sqrt{3}, -1+i\sqrt{3}, 1 ) $ of \cref{ex:whitehead}: i.e., it comes from the parameterisation
  \begin{displaymath}
    t \mapsto (\alpha,\beta,\lambda) = (2+2i\sqrt{3}, -1+i\sqrt{3}, t).
  \end{displaymath}
  In this example, step 2(b) of \cref{alg:stupid_picture} throws away all words which are in the parabolic subgroup $ \langle P,Q \rangle $ or its conjugates, because
  these words are parabolic at a test point that is known to lie in the interior of the space of compression body groups. The resulting picture is dark wherever there
  are `new' parabolics or elliptics arising in a particular group.
\end{ex}

\begin{figure}
  \centering
  \includegraphics[width=3.8in]{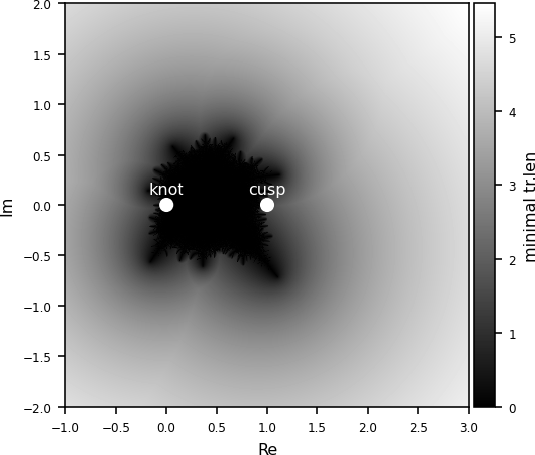}
  \caption{An affine slice through $ X((\Z \oplus \Z) * \Z) $ spanned by the Whitehead cusp and link groups.\label{fig:whitehead_cusp}}
\end{figure}
\begin{ex}[Whitehead link slice]\label{ex:whitehead_slice}
  In \cref{fig:whitehead_cusp}, we draw the affine slice of $ X((\Z \oplus \Z) * \Z) $ passing through the two groups produced in \cref{ex:whitehead}, the Whitehead link group and the corresponding
  cusp group. The slice is parameterised as a subset of the parameter space of groups $ G(\alpha,\beta,\lambda) $ by
  \begin{displaymath}
    t \mapsto (\alpha,\beta,\lambda) = \left((-2-2i) + (3+2i + i\sqrt{3})t,\, (1-i) - (3-i - 2i\sqrt{3})t,\, 1\right)\!,
  \end{displaymath}
  so $ t = 0 $ is the link group and $ t = 1 $ is the cusp group. The map $ \rho $ in \cref{alg:stupid_picture} is the map which sends $ t $
  to the triple $ \{P,Q,M\} $ with the values for $ \alpha $, $\beta$, and $ \lambda$ substituted in. Since $ M $ is parabolic for all groups in
  this slice, the quasiconformal deformation space of the generic discrete point in the image is a relatively open subset of the boundary
  of $ (1;2)$-compression body space in the character variety. The test point chosen for the algorithm must be different to the one chosen
  in \cref{ex:horiz_slice}, since $ M $ should be parabolic in addition to $ P $ and $ Q $.
\end{ex}

\begin{figure}
  \centering
  \includegraphics[width=4in]{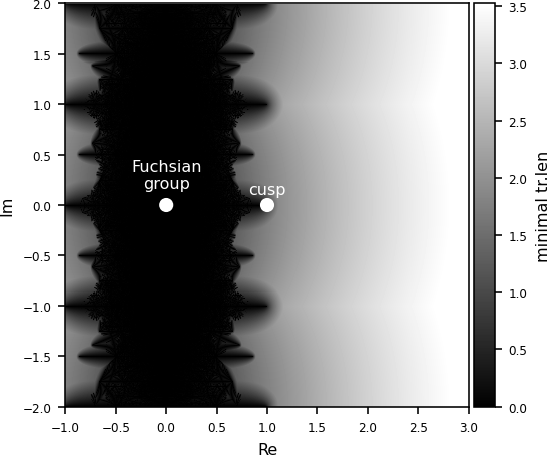}
  \caption{An affine slice through $ X((\Z \oplus \Z) * \Z) $ spanned by the Solomon's knot cusp group and a Fuchsian representation of the link group.\label{fig:solomon_cusp}}
\end{figure}
\begin{ex}[Solomon's knot slice]\label{ex:solomon_slice}
  Our last example in $X((\Z \oplus \Z) * \Z)$ is obtained from Solomon's knot, \cref{ex:solomon}. If we set up the system of equations
  \begin{displaymath}
    \tr M = 2, \tr P^{-1} M Q^{-1} = -2, MP^{-1} M Q^{-1} = \mathrm{Id}
  \end{displaymath}
  then all the solutions give groups conjugate to subgroups of $ \PSL(2,\R) $. The same is true for other choices of signs, showing that there is no discrete representation of the
  Solomon's knot group which is non-Fuchsian and providing an algebraic proof that the link is non-hyperbolic. We choose the non-faithful representation
  $ G( 2,2,1 ) $, and draw in \cref{fig:solomon_cusp} the affine slice through this point and the cusp point constructed in \cref{ex:solomon},
  \begin{displaymath}
    t \mapsto (\alpha,\beta,\lambda) = (2+2it, 2-2it, 1).
  \end{displaymath}
  Since both $ M $ and $ PQM^{-1} $ are parabolic for all groups in the image of this parameterisation, the large region of discrete groups in the figure
  is a subset of the boundary of $ (1;2)$-compression body space as in the previous example.
\end{ex}

As a final example, we return to the theme of reflection groups arising from gadgets made of rods and joints in the plane.

\begin{figure}
  \begin{subfigure}[m]{.35\textwidth}
    \centering
    \includegraphics[width=\textwidth]{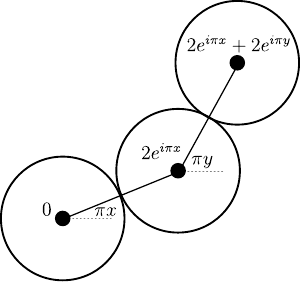}
    \caption{Parameters for the space.\label{fig:pendulum_diagram}}
  \end{subfigure}\hfill%
  \begin{subfigure}[m]{.6\textwidth}
    \centering
    \includegraphics[width=\textwidth]{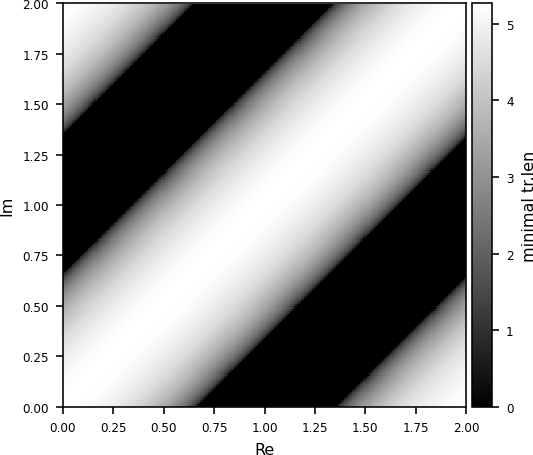}
    \caption{The space as a linear slice through $ X(F_3) $.\label{fig:pendulum}}
  \end{subfigure}
  \caption{The configuration space of a double pendulum.}
\end{figure}
\begin{ex}[The double pendulum]\label{ex:pendulum}
  Consider two rigid rods of length $2$, connected to form a double pendulum. Fixing one of the ends at $0$, the configuration space is parameterised
  by two angles $ \pi x $ and $ \pi y $ and we may draw circles of radius $1$ centred at the three vertices (\cref{fig:pendulum_diagram}). Computing with \cref{defn:reflection}, a
  generating set of the orientation-preserving half of the group generated by reflections in these circles (before normalising the determinant to $1$) is
  \begin{gather*}
    M_1 = \begin{bmatrix} -e^{i\pi x} & 2e^{2i\pi x}\\ -2 & 3e^{i\pi x} \end{bmatrix}\!,\;
    M_2 = \begin{bmatrix} -4 e^{i\pi x} - e^{i\pi y} & 8 e^{2i\pi x} + 8 e^{i\pi (x+y)} + 2 e^{2i\pi y}\\- 2 & 4 e^{i\pi x} + 3 e^{i\pi y} \end{bmatrix}\!,\,\text{and}\\
    M_3 = \begin{bmatrix} - 7 e^{i\pi (x+y)} - 4 e^{2i\pi y} - 4 e^{2i\pi x} & 2 e^{i\pi (2x+y)} + 2 e^{i\pi (x+2y)}\\ - 2 e^{i\pi x} - 2 e^{i\pi y} & e^{i\pi (x+y)} \end{bmatrix}\!.
  \end{gather*}
  Parameterising this family by the complex plane $ t = x+iy $ (with $ 0 \leq x \leq 2 $ and $ 0 \leq y \leq 2 $) and applying \cref{alg:stupid_picture} gives
  a picture of the phase space of the double pendulum, \cref{fig:pendulum}. The groups $ \langle M_1,M_2,M_3 \rangle $ generically lie in the boundary of genus $3$ Schottky space,
  and the lines in \cref{fig:pendulum} which appear straight are actually straight (not fractal), corresponding to the condition that the three circles become mutually tangent
  in a triangle.
\end{ex}

\section{Non-linear paths motivated by geometry}\label{sec:paths}
We began by studying groups generated by circles in the plane. We saw that if we take circles arising from lattice patterns then the resulting groups
can sometimes be remarkably close to knot groups in their presentation, and we described a way to formalise this in terms of unknotting tunnels. This led us
to consider linearly parameterised families of representations containing both knot groups and the corresponding cusp groups. However, there is no good reason to expect
these linear slices to respect any of the geometry involved; indeed, all of the trace algebra is highly non-linear. In this section we will restrict to a setting
that is quite well-understood, the setting of two-bridge links, where we can easily compute non-linear parameterisations that \emph{do} respect the geometry
that the knot groups and the cusp groups share.

\begin{figure}
  \centering
  \includegraphics[width=.7\textwidth]{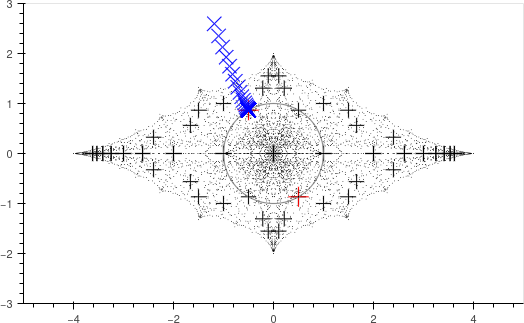}
  \caption{The shaded region is the exterior of $ \mc{R} $ and consists of indiscrete groups, $2$-bridge knot complement groups, and other Heckoid groups and their quotients. Knot complement
          groups with small denominator are indicated with $+$ signs. The path of points labelled $\times$ is the path along which the word representing the meridian of the figure 8 knot's
          unknotting tunnel is hyperbolic (the portion of the path in $ \mc{R} $) and elliptic of increasing holonomy angle (the portion from $\partial\mc{R}$ to
          the knot group, on the unit circle). \label{fig:fig8path}}
\end{figure}

\begin{thm}[{Riley \autocite{riley72}}]
  The holonomy group $ \pi $ of a $2$-bridge link $ \mf{b} $, indexed by a rational number $ p/q \in \Q $, has presentation $ \langle X, Y : W_{p/q} \rangle $
  where $ W_{p/q} $ is called the \df{Farey word} of slope $ p/q $. If $ \mf{b} $ is hyperbolic then there is a discrete faithful representation $ \rho : \pi \to \PSL(2,\C) $
  with $ \H^3/\rho(\pi) = \Sph^3 \setminus \mf{b} $, given by
  \begin{displaymath}
    \rho(X) = \begin{bmatrix} 1 & 1 \\ 0 & 1 \end{bmatrix} \;\text{and}\; \rho(Y) = \begin{bmatrix} 1 & 0 \\ z(p/q) & 1 \end{bmatrix}
  \end{displaymath}
  where $ z(p/q) $ is a solution in $ z $ to the equation $ \tr W_{p/q}(z) = 2 $. If $ \mf{b} $ is non-hyperbolic then there is still a discrete representation of this form, but
  $ z(p/q) \in \R $ and the image is Fuchsian. \qed
\end{thm}
Riley used this theorem (his version was phrased in terms of a factorisation of the Farey word rather than the full word itself) to conduct a computer study
of the representation space of $2$-bridge link groups, and produced a famous picture showing the locations of all $ z(p/q) \in \C $; a high resolution scan
of this image appears in Series \autocite{series19}. We discuss several aspects of this picture and the related theory in our expository article
with Martin and Schillewaert \autocite{ems22M}, and Lee and Sakuma carefully describe the structure of the discrete groups which appear in Riley's picture in their research announcement \autocite{lee12}.
The representation $ \rho : \pi \to \PSL(2,\C) $ lies in the character variety $ X(F_2) $; more precisely it lies inside the subvariety of representations
where the two generators are parabolic.
\begin{defn}[Riley slice]\label{defn:rileyslice}
  Consider the set of representations $ \rho : \langle X \rangle * \langle Y \rangle \to \PSL(2,\C) $ such that $ \rho(X) $ and $ \rho(Y) $ are parabolic. This
  set is a subvariety of $ X(F_2) $, and may be parameterised by a single complex number $ z $; we will choose the normalisation
  \begin{displaymath}
    \Gamma_z =  \left\langle \begin{bmatrix} 1 & 1 \\ 0 & 1 \end{bmatrix}\!,\, \begin{bmatrix} 1 & 0 \\ z & 1 \end{bmatrix} \right\rangle\!.
  \end{displaymath}
  When $ \abs{z} \gg 0 $, the group $ \Gamma_z $ is discrete
  and $ \Omega(\Gamma_z)/\Gamma_z $ is a four-punctured sphere; $ \H^3/\Gamma_z $ is a $3$-ball with two  ideal arcs drilled out. The set of all $ z $ such
  that $ \Gamma_z $ is of this form is called the \df{Riley slice}, $ \mc{R} $ \autocite{keen94}. It is an subset of $ \C $ which is conformally equivalent to a punctured disc
  (the puncture is at $ \infty $). The Riley slice is the exterior of the shaded `diamond' in \cref{fig:fig8path}.
\end{defn}

The reason we are spending so long on these definitions is that the Riley slice and the family of groups $ \Gamma_z $ is the most well-developed setting for a
theory connecting unknotting tunnels with knot representations and quasiconformal deformation spaces. A programme of study initiated by Sakuma \autocite{sakuma98}
and motivated in part by earlier work of Sakuma and Weeks \autocite{sakuma95} and the work of Riley discussed above conjectures that for each $ p/q $ there is
a smooth path through $ \C $, parameterised by $ \theta \in (0,2\pi) $, so that the point at $ \theta = 2\pi $ is
the $ p/q $ $2$-bridge knot, the point at $ \theta = 0 $ is the maximal cusp on $ \partial \mc{R} $ corresponding to pinching the four-punctured sphere along the curve with homology $ (p,q) $,
and the point at $ \theta \in (0,2\pi) $ corresponds to the cone manifold obtained by taking the complement of the $ p/q $ knot and replacing the upper unknotting
tunnel with a singular locus of angle $ \theta $. The curve is given by taking the Farey word $ W_{p/q}(z) $, which is the word
representing a loop around the unknotting tunnel (playing the same role as the dual curves we saw in the previous section), and choosing one of the branches
of $ \{ z \in \C : \tr W_{p/q}(z) \in (-2,2) \} $. The real trace locus of these words is already well-understood in the Riley slice interior \autocite{keen94,ems22M},
and the correct branch is known: it is the branch $ \mc{P}_{p/q} $ of $ \tr W_{p/q}^{-1}(\R) $ with the property that as $ x \to -\infty $, $ \Arg \tr W_{p/q}^{-1}(x) \to \pi p/q $.

Significant steps have been taken over the past two decades to realise this programme. The set of all discrete $ \Gamma_z $ has been fully classified \autocite{aimi2020classification,akiyoshi2020classification},
and from this classification we see that if $ z \in \mc{P}_{p/q} $ and $ \tr W_{p/q}(z) = 2\cos(\pi/n) $ then $ \Gamma_z $ is discrete and
uniformises the orbifold obtained by taking the knot complement of the $ p/q $ $2$-bridge link and replacing the upper unknotting tunnel with a singular arc of angle $ n $; and conversely,
any discrete group which is not in $ \overline{\mc{R}} $ is either the holonomy group of a $2$-bridge link complement with a singular unknotting tunnel (called a \df{Heckoid group}),
or the holonomy group of the quotient of such a manifold by a group of automorphisms.

\begin{figure}
  \begin{subfigure}[t]{.49\textwidth}
    \centering
    \includegraphics[width=\textwidth]{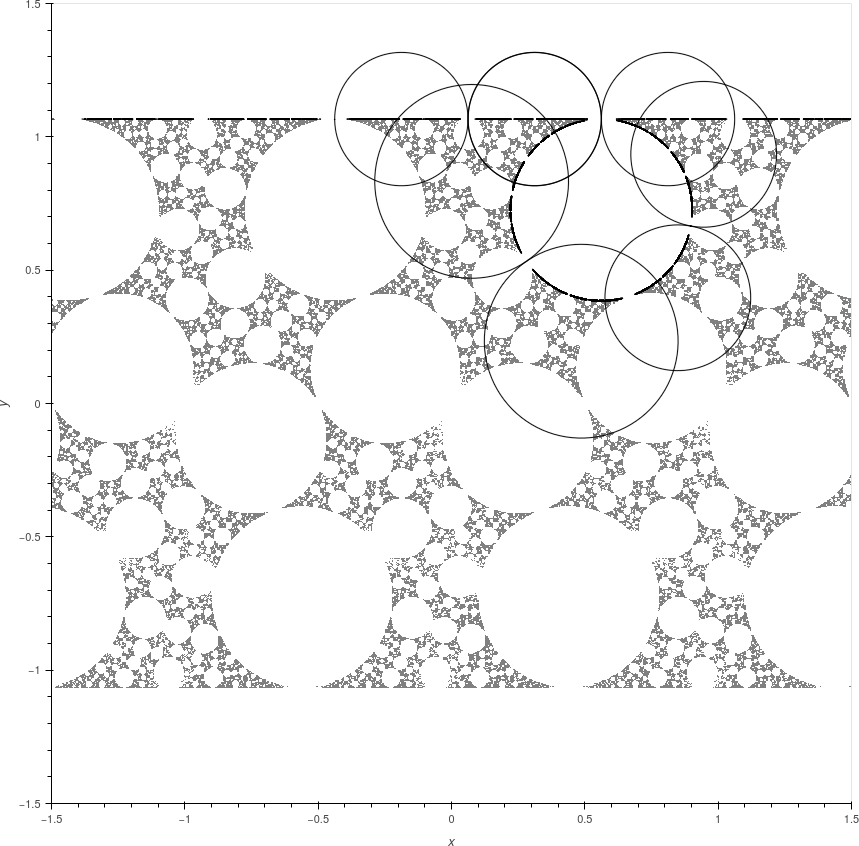}
    \caption{$z = -0.7733 + 1.4677i$, $ \theta = 0 $ (cusp)\label{fig:fig8limits_cusp}}
  \end{subfigure}\hfill%
  \begin{subfigure}[t]{.49\textwidth}
    \centering
    \includegraphics[width=\textwidth]{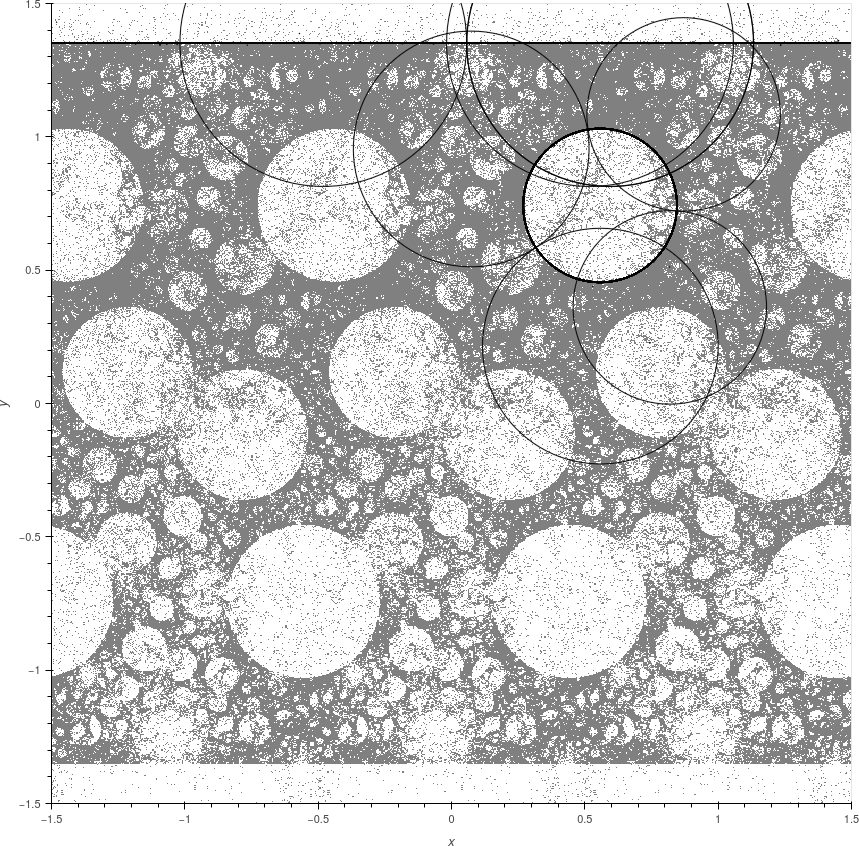}
    \caption{$z = -0.7151 + 1.3233i$, $ \theta = 1.045\pi $\label{fig:fig8limits_cone1}}
  \end{subfigure}\\[2ex]
  \begin{subfigure}[t]{.49\textwidth}
    \centering
    \includegraphics[width=\textwidth]{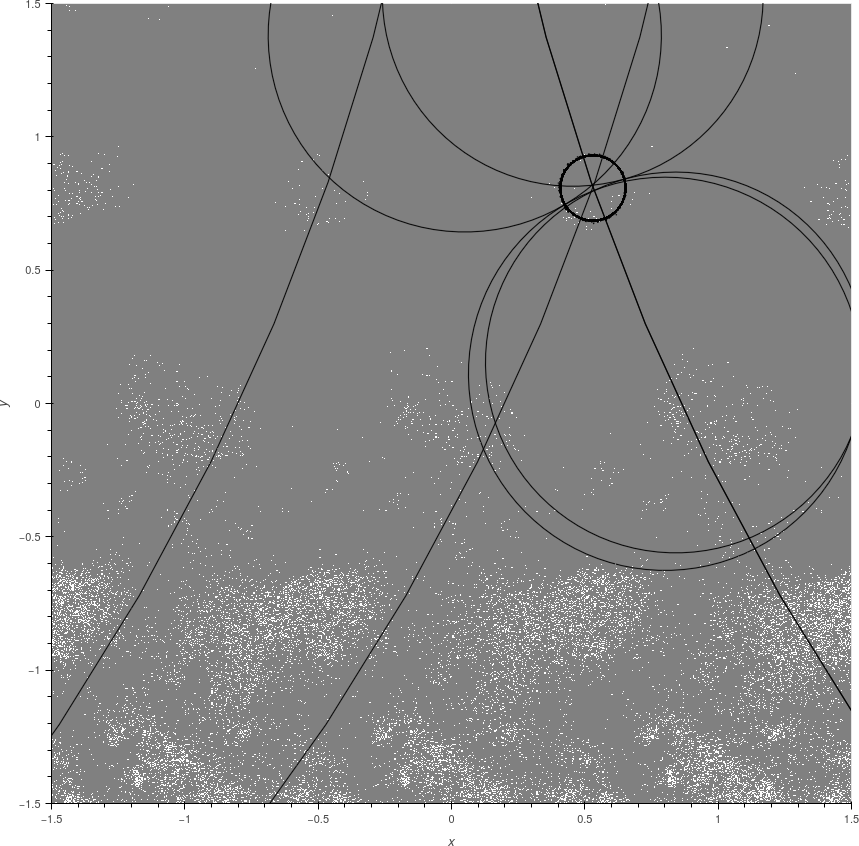}
    \caption{$z = -0.5768 + 1.0117i$, $ \theta = 1.786\pi $\label{fig:fig8limits_cone2}}
  \end{subfigure}\hfill%
  \begin{subfigure}[t]{.49\textwidth}
    \centering
    \includegraphics[width=\textwidth]{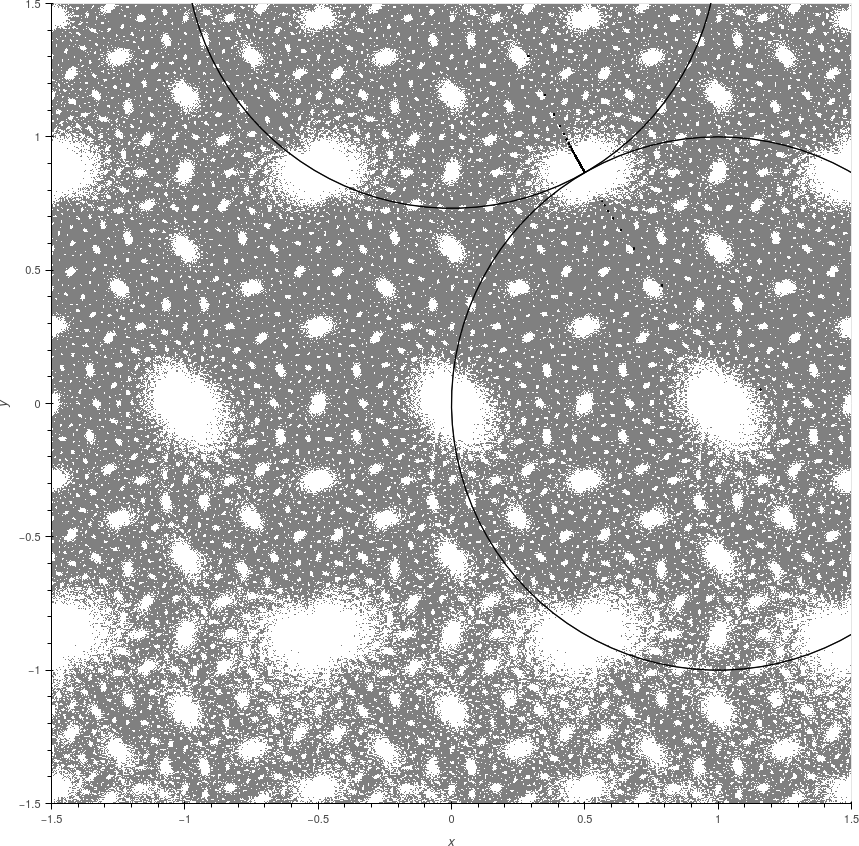}
    \caption{$z = -0.5 + 0.8660i$, $ \theta = 2\pi $ (knot group)\label{fig:fig8limits_knot}}
  \end{subfigure}
  \caption{Limit sets of conjectured cone manifold holonomy groups $ \Gamma_z $ with singular locus of angle $\theta $ along the figure eight knot's upper unknotting tunnel,
           showing limit points of embedded $S_{0,3}$ groups in black. The groups in \subref{fig:fig8limits_cone1} and \subref{fig:fig8limits_cone2}
           are indiscrete; since the elliptic elements of \subref{fig:fig8limits_cone1} have a smaller holonomy angle
           they require higher powers to be taken to fill the plane, which is why that picture is much less `dense'.\label{fig:fig8limits}}
\end{figure}

\begin{ex}[The figure eight knot]
  Consider the figure eight knot, which is the $2$-bridge link of slope $3/5$ and has Farey word $ W_{3/5} = XY^{-1} X^{-1} YXYX^{-1} Y^{-1} XY $. We will produce a family of limit sets passing
  down the conjectured path of cone manifolds. The algorithm that we use is elementary: we first find the point $ z_0 $ on some large circle around $0$, say radius $ R = 20 $, with the property
  that (i) the value of $ \tr W_{3/5} $ is real there, and (ii) it is the closest point on the circle to $ Re^{\pi i p/q} $ with this property. Set $ \tr W_{3/5} = x_0 $, and produce
  a list $ x_0, x_1,\ldots, x_n = -2 $ by taking $ n $ equally spaced points between $ 0 $ and $ \log (\abs{x_0}/2) $ and then applying the function $ -2\exp(\cdot) $ to all of these $n$ points;
  these numbers are traces which are (roughly) equally spaced with respect to the real translation length of $ W_{p/q} $. Now for each $ k $ from $ k = 1 $ to $ k = n $, find $ z_k $ by applying Newton's algorithm
  with initial guess $ z_{k-1} $ to solve the equation $ \tr W_{3/5}(z) = x_k $. We show this path superimposed on the Riley slice in \cref{fig:fig8path}, and we show how the limit set behaves
  along this path in \cref{fig:fig8limits}. Also shown in the figure are the limit sets of the two groups
  \begin{displaymath}
    \langle X, Y^{-1}X^{-1}YX Y X^{-1}Y^{-1}XY \rangle\;\text{and}\; \langle X Y^{-1}X^{-1}YXYX^{-1}, Y^{-1}XY \rangle
  \end{displaymath}
  which are the holonomy groups of the components of the convex core boundary of the figure eight knot cusp group and which degenerate to elementary rank one parabolic subgroups of the knot group.
\end{ex}

Similar constructions may be made for arbitrary tunnel number $1$ knots, using the theory of the previous section: from the knot we produce a cusp group on $ \partial \mc{S}_2 $,
with cusps corresponding to words representing curves bounding embedded discs in a handlebody in $ \Sph^3 $ which the knot union tunnel embeds into as a core $\theta$-curve. The systems
of equations corresponding to fixing the trace of the two meridian discs while continuously changing the trace of the dual curve to the unknotting tunnel are conjectured
to correspond to holonomy groups of hyperbolic cone manifolds that are homeomorphic to the knot complement and have a single singular arc, along the unknotting tunnel.

\sloppy
\printbibliography

@ARTICLE{accola66,
  AUTHOR = {Accola, Robert D.M.},
  DATE = {1966},
  DOI = {10.2307/2373196},
  JOURNALTITLE = {Amer. J. Math.},
  PAGES = {329--336},
  TITLE = {Invariant domains for Kleinian groups},
  VOLUME = {88},
}

@ARTICLE{aimi2020classification,
  AUTHOR = {Aimi, Shunsuke and Lee, Donghi and Sakai, Shunsuke and Sakuma, Makoto},
  DATE = {2020},
  EPRINT = {2001.11662},
  EPRINTCLASS = {math.GT},
  EPRINTTYPE = {arXiv},
  JOURNALTITLE = {Rend. Istit. Mat. Univ. Trieste},
  PAGES = {477--511},
  TITLE = {Classification of parabolic generating pairs of Kleinian groups with two parabolic generators},
  VOLUME = {52},
  DOI = {10.13137/2464-8728/30919}
}

@ARTICLE{akiyoshi2020classification,
  AUTHOR = {Akiyoshi, Hirotaka and Ohshika, Ken'ichi and Parker, John and Sakuma, Makoto and Yoshida, Han},
  DATE = {2021},
  DOI = {10.1090/tran/8246},
  EPRINT = {2001.09564},
  EPRINTCLASS = {math.GT},
  EPRINTTYPE = {arXiv},
  ISSUE = {3},
  JOURNALTITLE = {Trans. Amer. Math. Soc.},
  PAGES = {1765--1814},
  TITLE = {Classification of non-free Kleinian groups generated by two parabolic transformations},
  VOLUME = {374},
}

@BOOK{akiyoshi,
  AUTHOR = {Akiyoshi, Hirotaka and Sakuma, Makoto and Wada, Masaaki and Yamashita, Yasushi},
  PUBLISHER = {Springer},
  DATE = {2007},
  DOI = {10.1007/978-3-540-71807-9},
  ISBN = {9783540718079},
  NUMBER = {1909},
  SERIES = {Lecture Notes in Math.},
  TITLE = {Punctured torus groups and $2$-bridge knot groups I},
}

@BOOK{apanasov,
  AUTHOR = {Apanasov, Boris N.},
  PUBLISHER = {Kluwer Academic Publishers},
  DATE = {1991},
  ISBN = {0-7923-0216-8},
  NUMBER = {40},
  SERIES = {Math. Appl. (Soviet Ser.)},
  TITLE = {Discrete groups in space and uniformization problems},
}

@ARTICLE{apanasov80r,
  AUTHOR = {Apanasov, Boris N.},
  DATE = {1980},
  JOURNALTITLE = {Sibirsk. Mat. Zh.},
  NUMBER = {4},
  PAGES = {3--15},
  TITLE = {Kleinian groups, Teichmüller space and Mostow's rigidity theorem},
  VOLUME = {21},
  LANGUAGE = {Russian},
  URL = {https://www.mathnet.ru/php/archive.phtml?wshow=paper&jrnid=smj&paperid=3745&option_lang=eng}
}

@ARTICLE{apanasov80e,
  AUTHOR = {Apanasov, Boris N.},
  DATE = {1981},
  JOURNALTITLE = {Siberian Math. J.},
  NUMBER = {4},
  PAGES = {483--491},
  TITLE = {Kleinian groups, Teichmüller space and Mostow's rigidity theorem},
  VOLUME = {21},
  DOI = {10.1007/BF00995946}
}

@BOOK{bogatyrev,
  AUTHOR = {Bogatyrev, Andrei},
  PUBLISHER = {Springer},
  DATE = {2012},
  ISBN = {978-3-642-25633-2},
  NUMBER = {149},
  SERIES = {Springer Monogr. Math.},
  TITLE = {Extremal polynomials and Riemann surfaces},
}

@BOOK{bottazzini_gray,
  AUTHOR = {Bottazzini, Umberto and Gray, Jeremy},
  PUBLISHER = {Springer},
  DATE = {2013},
  DOI = {https://doi.org/10.1007/978-1-4614-5725-1},
  ISBN = {978-1-4614-5725-1},
  SERIES = {Sources Stud. Hist. Math. Phys. Sci.},
  SUBTITLE = {The rise of complex function theory},
  TITLE = {Hidden harmony---Geometric fantasies},
}

@ARTICLE{Brin13,
  AUTHOR = {Brin, Matthew G. and Jones, Gareth A. and Singerman, David},
  DATE = {2013},
  DOI = {10.1016/j.exmath.2013.01.002},
  JOURNALTITLE = {Expo. Math.},
  NUMBER = {2},
  PAGES = {99--103},
  TITLE = {Commentary on Robert Riley's article `A personal account of the discovery of hyperbolic structures on some knot complements'},
  VOLUME = {31},
}

@ARTICLE{brooks86,
  AUTHOR = {Brooks, Robert},
  DATE = {1986},
  DOI = {10.1007/BF01389263},
  JOURNALTITLE = {Invent. Math.},
  PAGES = {461--469},
  TITLE = {Circle packings and co-compact extensions of Kleinian groups},
  VOLUME = {86},
}

@MISC{bulatov13,
  AUTHOR = {Bulatov, Vladimir},
  URL = {https://gallery.bridgesmathart.org/exhibitions/2013-joint-mathematics-meetings/bulatov},
  DATE = {2012},
  HOWPUBLISHED = {Digital print},
  TITLE = {Klein-Fricke composition group tiling},
}

@ARTICLE{cho07,
  AUTHOR = {Cho, Sangbum and McCullough, Darryl},
  DATE = {2009},
  DOI = {10.2140/gt.2009.13.769},
  EPRINT = {math/0611921},
  EPRINTCLASS = {math.GT},
  EPRINTTYPE = {arXiv},
  JOURNALTITLE = {Geom. Topol.},
  PAGES = {769--815},
  TITLE = {The tree of knot tunnels},
  VOLUME = {13},
}

@BOOK{coxAPS,
  AUTHOR = {Cox, David A.},
  PUBLISHER = {American Mathematical Society},
  DATE = {2020},
  NUMBER = {134},
  SERIES = {CBMS Reg. Conf. Ser. Math.},
  TITLE = {Applications of polynomial systems},
}

@ARTICLE{dang19,
  AUTHOR = {Dang, Vinh and Purcell, Jessica S.},
  DATE = {2019},
  DOI = {10.1090/proc/14336},
  EPRINT = {1711.03693},
  EPRINTCLASS = {math.GT},
  EPRINTTYPE = {arXiv},
  JOURNALTITLE = {Proc. Amer. Math. Soc.},
  PAGES = {1351--1366},
  TITLE = {Cusp shape and tunnel number},
  VOLUME = {147},
}

@BOOK{dehn,
  AUTHOR = {Dehn, Max},
  TRANSLATOR = {Stillwell, John},
  PUBLISHER = {Springer-Verlag},
  DATE = {1987},
  ISBN = {0387964169},
  TITLE = {Papers on group theory and topology},
}

@ARTICLE{diaz19,
  AUTHOR = {Díaz, Juan Pablo and Hinojosa, Gabriela and Mendoza, Martha and Verjovsky, Alberto},
  DATE = {2019},
  DOI = {10.1080/17513472.2018.1506615},
  JOURNALTITLE = {J. Math. Arts},
  NUMBER = {3},
  PAGES = {230--242},
  TITLE = {Dynamically defined wild knots and Othoniel's \emph{My Way}},
  VOLUME = {13},
}

@MISC{elzenaar24c,
  AUTHOR = {Elzenaar, Alex},
  DATE = {2024},
  EPRINT = {2411.17940},
  EPRINTCLASS = {math.GT},
  EPRINTTYPE = {arXiv},
  TITLE = {Changing topological type of compression bodies through cone manifolds},
}

@ARTICLE{ems21,
  AUTHOR = {Elzenaar, Alex and Martin, Gaven J. and Schillewaert, Jeroen},
  DATE = {2023},
  DOI = {10.1016/j.exmath.2022.12.002},
  EPRINT = {2111.03230},
  EPRINTCLASS = {math.GT},
  EPRINTTYPE = {arXiv},
  ISSUE = {1},
  JOURNALTITLE = {Expo. Math.},
  PAGES = {20--54},
  TITLE = {Approximations of the Riley slice},
  VOLUME = {41},
}

@INBOOK{ems22M,
  AUTHOR = {Elzenaar, Alex and Martin, Gaven J. and Schillewaert, Jeroen},
  EDITOR = {Wood, David R. and de Gier, Jan and Praeger, Cheryl E.},
  PUBLISHER = {Springer},
  BOOKTITLE = {2021-22 MATRIX annals},
  DATE = {2024},
  DOI = {10.1007/978-3-031-47417-0_2},
  EPRINT = {2204.11422},
  EPRINTCLASS = {math.GT},
  EPRINTTYPE = {arXiv},
  ISBN = {978-3-031-47417-0},
  PAGES = {31--74},
  TITLE = {Concrete one complex dimensional moduli spaces of hyperbolic manifolds and orbifolds},
}

@MISC{fairchild24,
  AUTHOR = {Fairchild, Samantha and Ortiz, Ángel David Ríos},
  DATE = {2024},
  EPRINT = {2401.10801},
  EPRINTCLASS = {math.GT},
  EPRINTTYPE = {arXiv},
  TITLE = {Crossing the transcendental divide: from Schottky groups to algebraic curves},
}

@BOOK{fricke_klein65t,
  AUTHOR = {Fricke, Robert and Klein, Felix},
  TRANSLATOR = {DuPre, Arthur M.},
  PUBLISHER = {Higher Education Press},
  DATE = {2017},
  ISBN = {9787040478402},
  NUMBER = {3},
  RELATED = {fricke_klein65},
  RELATEDTYPE = {translationof},
  SERIES = {CTM. Class. Top. Math.},
  TITLE = {Lectures on the theory of automorphic functions},
}

@BOOK{fricke_klein65,
  AUTHOR = {Fricke, Robert and Klein, Felix},
  LOCATION = {Stuttgart},
  PUBLISHER = {B. G. Teubner},
  DATE = {1897},
  RELATED = {c4b8818a22d85ac62c36b7ec876f8408},
  RELATEDTYPE = {translatedas},
  TITLE = {Vorlesungen über die Theorie der automorphen Funktionen},
}

@BOOK{grayP,
  AUTHOR = {Gray, Jeremy},
  PUBLISHER = {Princeton University Press},
  DATE = {2013},
  ISBN = {978-0-691-15271-4},
  SUBTITLE = {A scientific biography},
  TITLE = {Henri Poincaré},
}

@BOOK{grunbaum,
  AUTHOR = {Grünbaum, Branko and Shephard, G. C.},
  PUBLISHER = {Dover Publications},
  DATE = {2016},
  EDITION = {2},
  ISBN = {9780486469812},
  TITLE = {Tilings \& patterns},
}

@INBOOK{hk03,
  AUTHOR = {Hodgson, Craig D. and Kerckhoff, Steven P.},
  PUBLISHER = {Cambridge University Press},
  BOOKTITLE = {Kleinian groups and hyperbolic $3$-manifolds},
  DATE = {2003},
  NUMBER = {299},
  PAGES = {41--74},
  SERIES = {London Math. Soc. Lecture Note Ser.},
  TITLE = {Harmonic deformations of hyperbolic $3$-manifolds},
}

@ARTICLE{hk98,
  AUTHOR = {Hodgson, Craig D. and Kerckhoff, Steven P.},
  URL = {http://projecteuclid.org/euclid.jdg/1214460606},
  DATE = {1998},
  JOURNALTITLE = {J. Differential Geom.},
  PAGES = {1--59},
  TITLE = {Rigidity of hyperbolic cone-manifolds and hyperbolic Dehn surgery},
  VOLUME = {48},
}

@ARTICLE{horowitz72,
  AUTHOR = {Horowitz, Robert D.},
  DATE = {1972},
  DOI = {10.1002/cpa.3160250602},
  JOURNALTITLE = {Comm. Pure Appl. Math.},
  PAGES = {635--649},
  TITLE = {Characters of free groups represented in the two-dimensional special linear group},
  VOLUME = {25},
}

@ARTICLE{hyde03,
  AUTHOR = {Hyde, Stephen T. and Larsson, Ann-Kristin and Matteo, Tiziana Di and Ramsden, Stuart and Robins, Vanessa},
  DATE = {2003},
  DOI = {10.1071/CH03191},
  JOURNALTITLE = {Aust. J. Chem.},
  PAGES = {981--1000},
  TITLE = {Meditation on an engraving of Fricke and Klein (The modular group and geometrical chemistry)},
  VOLUME = {56},
}

@MISC{ichikawa25,
  AUTHOR = {Ichikawa, Takashi and Kodama, Yuji},
  DATE = {2025},
  EPRINT = {2507.20296},
  EPRINTCLASS = {nlin.SI},
  EPRINTTYPE = {arXiv},
  TITLE = {KP solitons and the Schottky uniformization},
}

@ARTICLE{ishihara11,
  AUTHOR = {Ishahara, Kai},
  DATE = {2011},
  DOI = {10.2140/agt.2011.11.2167},
  JOURNALTITLE = {Algebr. Geom. Topol.},
  PAGES = {2167--2190},
  TITLE = {An algorithm for finding parameters of tunnels},
  VOLUME = {11},
}

@BOOK{kapovich,
  AUTHOR = {Kapovich, Michael},
  PUBLISHER = {Birkhaüser},
  DATE = {2001},
  DOI = {10.1007/978-0-8176-4913-5},
  ISBN = {978-0-8176-4913-5},
  NUMBER = {183},
  SERIES = {Progr. Math.},
  TITLE = {Hyperbolic manifolds and discrete groups},
}

@ARTICLE{kapovich23,
  AUTHOR = {Kapovich, Michael and Kontorovich, Alex},
  DATE = {2023},
  DOI = {10.1515/crelle-2023-0004},
  EPRINT = {2104.13838},
  EPRINTCLASS = {math.NT},
  EPRINTTYPE = {arXiv},
  JOURNALTITLE = {J. Reine Angew. Math.},
  PAGES = {105--142},
  TITLE = {On superintegral Kleinian sphere packings, bugs, and arithmetic groups},
  VOLUME = {798},
}

@ARTICLE{keen91,
  AUTHOR = {Keen, Linda and Maskit, Bernard and Series, Caroline},
  DATE = {1993},
  DOI = {10.1515/crll.1993.436.209},
  EPRINT = {math/9201299},
  EPRINTCLASS = {math.DG},
  EPRINTTYPE = {arXiv},
  JOURNALTITLE = {J. Reine Angew. Math.},
  PAGES = {209--219},
  TITLE = {Geometric finiteness and uniqueness for Kleinian groups with circle packing limit sets},
  VOLUME = {436},
}

@ARTICLE{keen94,
  AUTHOR = {Keen, Linda and Series, Caroline},
  DATE = {1994},
  DOI = {10.1112/plms/s3-69.1.72},
  JOURNALTITLE = {Proc. London Math. Soc. (3)},
  NUMBER = {1},
  PAGES = {72--90},
  TITLE = {The Riley slice of Schottky space},
  VOLUME = {69},
}

@BOOK{fricke_klein66t,
  AUTHOR = {Klein, Felix and Fricke, Robert},
  TRANSLATOR = {DuPre, Arthur M.},
  PUBLISHER = {Higher Education Press},
  DATE = {2017},
  ISBN = {9787040478372},
  NUMBER = {3},
  RELATED = {fricke_klein66},
  RELATEDTYPE = {translationof},
  SERIES = {CTM. Class. Top. Math.},
  TITLE = {Lectures on the theory of elliptic modular functions},
}

@BOOK{fricke_klein66,
  AUTHOR = {Klein, Felix and Fricke, Robert},
  LOCATION = {Leipzig},
  PUBLISHER = {B. G. Teubner},
  DATE = {1890},
  RELATED = {aa5baba684cb410bba1a59f63cc47c84},
  RELATEDTYPE = {translatedas},
  TITLE = {Vorlesungen über die Theorie der elliptischen Modulfunctionen},
}

@ARTICLE{komori06,
  AUTHOR = {Komori, Yohei and Sugawa, Toshiyuki and Wada, Masaaki and Yamashita, Yasushi},
  DATE = {2006},
  DOI = {10.1080/10586458.2006.10128951},
  JOURNALTITLE = {Exp. Math.},
  NUMBER = {1},
  PAGES = {51--60},
  TITLE = {Drawing Bers embeddings of the Teichmüller space of once-punctured tori},
  VOLUME = {15},
}

@BOOK{kag,
  AUTHOR = {Krushkaĺ, S. L. and Apanasov, B. N. and Gusevskiĭ, N. A.},
  EDITOR = {Maskit, Bernard},
  TRANSLATOR = {McFaden, H. H.},
  PUBLISHER = {American Mathematical Society},
  DATE = {1986},
  NUMBER = {62},
  SERIES = {Transl. Math. Monogr.},
  TITLE = {Kleinian groups and uniformization in examples and problems},
}

@ARTICLE{lackenbypurcell13,
  AUTHOR = {Lackenby, Marc and Purcell, Jessica S.},
  DATE = {2014},
  DOI = {10.1080/10586458.2013.870503},
  EPRINT = {1302.3652},
  EPRINTCLASS = {math.GT},
  EPRINTTYPE = {arXiv},
  JOURNALTITLE = {Exp. Math.},
  PAGES = {218--240},
  TITLE = {Geodesics and compression bodies},
  VOLUME = {23},
}

@ARTICLE{lee12,
  AUTHOR = {Lee, Donghi and Sakuma, Makoto},
  DATE = {2012},
  DOI = {10.3934/era.2012.19.97},
  EPRINT = {1206.4258},
  EPRINTCLASS = {math.GT},
  EPRINTTYPE = {arXiv},
  JOURNALTITLE = {Electron. Res. Announc. Math. Sci.},
  PAGES = {97--111},
  TITLE = {Simple loops on $2$-bridge spheres in Heckoid orbifolds for $2$-bridge links},
  VOLUME = {19},
}

@BOOK{magnus,
  AUTHOR = {Magnus, Wilhelm},
  PUBLISHER = {Academic Press},
  DATE = {1974},
  ISBN = {0124654509},
  NUMBER = {61},
  SERIES = {Pure Appl. Math.},
  TITLE = {Noneuclidean tesselations and their groups},
}

@BOOK{mandelbrotCW,
  AUTHOR = {Mandelbrot, Benoit B.},
  PUBLISHER = {Springer},
  DATE = {2004},
  DOI = {10.1007/978-1-4757-4017-2},
  ISBN = {978-1-4757-4017-2},
  SUBTITLE = {The Mandelbrot set and beyond},
  TITLE = {Fractals and chaos},
}

@BOOK{marden,
  AUTHOR = {Marden, Albert},
  PUBLISHER = {Cambridge University Press},
  DATE = {2016},
  EDITION = {2},
  ISBN = {9781107116740},
  NOTE = {First edition was published under the title ``Outer circles''},
  SUBTITLE = {An introduction in 2 and 3 dimensions},
  TITLE = {Hyperbolic manifolds},
}

@BOOK{maskit,
  AUTHOR = {Maskit, Bernard},
  PUBLISHER = {Springer-Verlag},
  DATE = {1987},
  DOI = {10.1007/978-3-642-61590-0},
  ISBN = {978-3-642-61590-0},
  NUMBER = {287},
  SERIES = {Grundlehren Math. Wiss.},
  TITLE = {Kleinian groups},
}

@BOOK{matsuzakitaniguchi,
  AUTHOR = {Matsuzaki, Katsuhiko and Taniguchi, Masahiko},
  PUBLISHER = {Oxford University Press},
  DATE = {1998},
  ISBN = {0198500629},
  TITLE = {Hyperbolic manifolds and Kleinian groups},
}

@BOOK{mccolloughmiller,
  AUTHOR = {McCollough, Darryl and Miller, Andy},
  PUBLISHER = {American Mathematical Society},
  DATE = {1986},
  NUMBER = {344},
  SERIES = {Mem. Amer. Math. Soc.},
  TITLE = {Homeomorphisms of $3$-manifolds with compressible boundary},
}

@BOOK{mcmullenRFC,
  AUTHOR = {McMullen, Curtis T.},
  PUBLISHER = {Princeton University Press},
  DATE = {1996},
  NUMBER = {142},
  SERIES = {Ann. of Math. Stud.},
  TITLE = {Renormalization and $3$-manifolds which fiber over the circle},
}

@ARTICLE{yokota96,
  AUTHOR = {Morimoto, Kanji and Sakuma, Makoto and Yokota, Yoshiyuki},
  DATE = {1996},
  DOI = {10.2969/jmsj/04840667},
  JOURNALTITLE = {J. Math. Soc. Japan},
  NUMBER = {4},
  PAGES = {667--688},
  TITLE = {Identifying tunnel number one knots},
  VOLUME = {48},
}

@BOOK{indras_pearls,
  AUTHOR = {Mumford, David and Series, Caroline and Wright, David},
  PUBLISHER = {Cambridge University Press},
  DATE = {2002},
  ISBN = {0521352533},
  SUBTITLE = {The vision of {F}elix {K}lein},
  TITLE = {Indra's pearls},
}

@ARTICLE{nakamura21,
  AUTHOR = {Nakamura, Kento},
  DATE = {2021},
  DOI = {10.1080/17513472.2021.1943998},
  JOURNALTITLE = {J. Math. Arts},
  NUMBER = {2},
  PAGES = {106--136},
  TITLE = {Iterated inversion system: an algorithm for efficiently visualizing Kleinian groups and extending the possibilities of fractal art},
  VOLUME = {15},
}

@INBOOK{nimershiem94,
  AUTHOR = {Nimershiem, Barbara E.},
  EDITOR = {Johannson, Klaus},
  PUBLISHER = {International Press},
  BOOKTITLE = {Low-dimensional topology (Knoxville, TN, 1992)},
  DATE = {1994},
  ISBN = {1-57146-018-7},
  NUMBER = {III},
  PAGES = {133--142},
  SERIES = {Conf. Proc. Lecture Notes Geom. Topology},
  TITLE = {Isometry classes of flat $2$-tori appearing as cusps of hyperbolic $3$-manifolds are dense in the moduli space of the torus},
}

@ARTICLE{page15,
  AUTHOR = {Page, Aurel},
  URL = {https://www.jstor.org/stable/24489209},
  DATE = {2015},
  JOURNALTITLE = {Math. Comp.},
  NUMBER = {295},
  PAGES = {2361--2390},
  TITLE = {Computing arithmetic Kleinian groups},
  VOLUME = {84},
}

@BOOK{poincare,
  AUTHOR = {Poincaré, Henri},
  TRANSLATOR = {Stillwell, John},
  PUBLISHER = {Springer},
  DATE = {1985},
  TITLE = {Papers on Fuchsian functions},
}

@INBOOK{rigby95,
  AUTHOR = {Rigby, J.F.},
  EDITOR = {Emmer, Michele},
  PUBLISHER = {Leonardo Books},
  BOOKSUBTITLE = {Art and mathematics},
  BOOKTITLE = {The visual mind},
  CHAPTER = {26},
  DATE = {1993},
  ISBN = {9780262050487},
  PAGES = {177--186},
  TITLE = {Compound tilings and perfect colourings},
}

@ARTICLE{Riley13,
  AUTHOR = {Riley, Robert},
  DATE = {2013},
  DOI = {10.1016/j.exmath.2013.01.003},
  EPRINT = {1301.4601},
  EPRINTTYPE = {arXiv},
  JOURNALTITLE = {Expo. Math.},
  NUMBER = {2},
  PAGES = {104--115},
  TITLE = {A personal account of the discovery of hyperbolic structures on some knot complements},
  VOLUME = {31},
}

@ARTICLE{riley75b,
  AUTHOR = {Riley, Robert},
  DATE = {1975},
  DOI = {10.1017/s0305004100051094},
  JOURNALTITLE = {Math. Proc. Cambridge Philos. Soc.},
  PAGES = {281--288},
  TITLE = {A quadratic parabolic group},
  VOLUME = {77},
}

@INBOOK{riley79,
  AUTHOR = {Riley, Robert},
  EDITOR = {Fenn, Roger},
  PUBLISHER = {Springer-Verlag},
  BOOKTITLE = {Topology of low-dimensional manifolds---Proceedings, Sussex 1977},
  DATE = {1979},
  DOI = {10.1007/BFb0063194},
  ISBN = {3540095063},
  NUMBER = {722},
  PAGES = {99--133},
  SERIES = {Lecture Notes in Math.},
  TITLE = {An elliptical path from parabolic representations to hyperbolic structures},
}

@ARTICLE{riley83,
  AUTHOR = {Riley, Robert},
  DATE = {1983},
  DOI = {10.2307/2007537},
  JOURNALTITLE = {Math. Comp.},
  NUMBER = {162},
  PAGES = {607--632},
  TITLE = {Applications of a computer implementation of Poincaré's theorem on fundamental polyhedra},
  VOLUME = {40},
}

@ARTICLE{riley75c,
  AUTHOR = {Riley, Robert},
  DATE = {1975},
  DOI = {10.1112/S0025579300005982},
  JOURNALTITLE = {Mathematika},
  PAGES = {141--150},
  TITLE = {Discrete parabolic representations of link groups},
  VOLUME = {22},
}

@ARTICLE{riley72,
  AUTHOR = {Riley, Robert},
  DATE = {1972},
  DOI = {10.1112/plms/s3-24.2.217},
  JOURNALTITLE = {Proc. London Math. Soc. (3)},
  PAGES = {217--242},
  TITLE = {Parabolic representations of knot groups, I},
  VOLUME = {24},
}

@ARTICLE{riley75,
  AUTHOR = {Riley, Robert},
  DATE = {1975},
  DOI = {10.1112/plms/s3-31.4.495},
  ISSUE = {4},
  JOURNALTITLE = {Proc. London Math. Soc. (3)},
  PAGES = {495--512},
  TITLE = {Parabolic representations of knot groups, II},
  VOLUME = {31},
}

@INPROCEEDINGS{riley82,
  AUTHOR = {Riley, Robert},
  EDITOR = {Brown, R. and Thickstun, T.L.},
  PUBLISHER = {Cambridge University Press},
  BOOKTITLE = {Low-dimensional topology},
  DATE = {1982},
  DOI = {10.1017/CBO9780511758935.009},
  EVENTDATE = {1979},
  EVENTTITLE = {Conference on Topology in Low Dimension},
  ISBN = {0-521-28146-6},
  NUMBER = {48},
  PAGES = {81--151},
  SERIES = {London Math. Soc. Lecture Note Ser.},
  TITLE = {Seven excellent knots},
  VENUE = {University College of North Wales, Bangor},
}

@INBOOK{rojas98,
  AUTHOR = {Rojas, Raúl},
  EDITOR = {Pickover, Clifford A.},
  PUBLISHER = {Elsevier},
  BOOKSUBTITLE = {Ten year compilation of advanced research},
  BOOKTITLE = {Chaos and fractals: A computer graphical journey},
  DATE = {1998},
  ISBN = {0444500022},
  PAGES = {225--234},
  TITLE = {A tutorial on efficient computer graphic representations of the Mandlebrot set},
}

@ARTICLE{sakuma98,
  AUTHOR = {Sakuma, Makoto},
  DATE = {1998-04},
  DOI = {10.1016/s0960-0779(97)00101-x},
  JOURNALTITLE = {Chaos Solitons Fractals},
  NUMBER = {4-5},
  PAGES = {739--748},
  TITLE = {The topology, geometry and algebra of unknotting tunnels},
  VOLUME = {9},
}

@ARTICLE{sakuma95,
  AUTHOR = {Sakuma, Makoto and Weeks, Jeffrey},
  DATE = {1995},
  DOI = {10.4099/math1924.21.393},
  JOURNALTITLE = {Japan. J. Math. (N.S.)},
  NUMBER = {2},
  PAGES = {393--439},
  TITLE = {Examples of canonical decompositions of hyperbolic link complements},
  VOLUME = {21},
}

@MISC{series19,
  AUTHOR = {Series, Caroline},
  URL = {https://warwick.ac.uk/fac/sci/maths/people/staff/caroline_series/preslecture.pdf},
  DATE = {2019},
  HOWPUBLISHED = {Beamer slides},
  TITLE = {All about the Riley slice: LMS Presidential Lecture 2019},
}

@ARTICLE{series17,
  AUTHOR = {Series, Caroline and Tan, Ser and Yamashita, Yasushi},
  DATE = {2017},
  DOI = {10.2140/agt.2017.17.2239},
  EPRINT = {1409.6863},
  EPRINTTYPE = {arXiv},
  JOURNALTITLE = {Algebr. Geom. Topol.},
  NUMBER = {4},
  PAGES = {2239--2282},
  TITLE = {The diagonal slice of Schottky space},
  VOLUME = {17},
}

@MISC{stange24,
  AUTHOR = {Stange, Katherine E.},
  DATE = {2024},
  EPRINT = {2412.02050},
  EPRINTCLASS = {math.NT},
  EPRINTTYPE = {arXiv},
  TITLE = {An illustrated introduction to the arithmetic of Apollonian circle packings, continued fractions, and other thin orbits},
}

@MISC{stier17,
  AUTHOR = {Stier, Torsten},
  URL = {https://artscience-node.com/wp-content/uploads/2017/06/Lange-Nacht-2017_Stier_with-heading.pdf},
  DATE = {2017},
  HOWPUBLISHED = {3D printed sculpture},
  TITLE = {42},
}

@ARTICLE{wada06,
  AUTHOR = {Wada, Masaaki},
  URL = {http://projecteuclid.org/euclid.em/1150476904},
  DATE = {2006},
  JOURNALTITLE = {Exp. Math.},
  PAGES = {61--66},
  TITLE = {OPTi's algorithm for discreteness determination},
  VOLUME = {15},
}

@ARTICLE{wielenberg78,
  AUTHOR = {Wielenberg, Norbert},
  DATE = {1978},
  DOI = {10.1017/s0305004100055250},
  JOURNALTITLE = {Math. Proc. Cambridge Philos. Soc.},
  PAGES = {427--436},
  TITLE = {The structure of certain subgroups of the Picard group},
  VOLUME = {84},
}

@INBOOK{wright05,
  AUTHOR = {Wright, David},
  EDITOR = {Minsky, Yair N. and Sakuma, Makoto and Series, Caroline},
  PUBLISHER = {Cambridge University Press},
  BOOKTITLE = {Spaces of Kleinian groups},
  DATE = {2006},
  DOI = {10.1017/cbo9781139106993.016},
  ISBN = {9780521617970},
  NUMBER = {329},
  PAGES = {301--336},
  SERIES = {London Math. Soc. Lecture Note Ser.},
  TITLE = {Searching for the cusp},
}

@INBOOK{yamashita12,
  AUTHOR = {Yamashita, Yasushi},
  EDITOR = {Goldman, William and Series, Caroline and Tan, Ser Peow},
  PUBLISHER = {World Scientific},
  BOOKTITLE = {Geometry, topology and dynamics of character varieties},
  DATE = {2012},
  DOI = {10.1142/9789814401364_0005},
  NUMBER = {23},
  PAGES = {159--190},
  SERIES = {Lect. Notes Ser. Inst. Math. Sci. Natl. Univ. Singap.},
  TITLE = {Creating software for visualizing Kleinian groups},
}

\end{document}